

%
\magnification1200
\pretolerance=100
\tolerance=200
\hbadness=1000
\vbadness=1000
\linepenalty=10
\hyphenpenalty=50
\exhyphenpenalty=50
\binoppenalty=700
\relpenalty=500
\clubpenalty=5000
\widowpenalty=5000
\displaywidowpenalty=50
\brokenpenalty=100
\predisplaypenalty=7000
\postdisplaypenalty=0
\interlinepenalty=10
\doublehyphendemerits=10000
\finalhyphendemerits=10000
\adjdemerits=160000
\uchyph=1
\delimiterfactor=901
\hfuzz=0.1pt
\vfuzz=0.1pt
\overfullrule=5pt
\hsize=146 true mm
\vsize=8.9 true in
\maxdepth=4pt
\delimitershortfall=.5pt
\nulldelimiterspace=1.2pt
\scriptspace=.5pt
\normallineskiplimit=.5pt
\mathsurround=0pt
\parindent=20pt
\catcode`\_=11
\catcode`\_=8
\normalbaselineskip=12pt
\normallineskip=1pt plus .5 pt minus .5 pt
\parskip=6pt plus 3pt minus 3pt
\abovedisplayskip = 12pt plus 5pt minus 5pt
\abovedisplayshortskip = 1pt plus 4pt
\belowdisplayskip = 12pt plus 5pt minus 5pt
\belowdisplayshortskip = 7pt plus 5pt
\normalbaselines
\smallskipamount=\parskip
 \medskipamount=2\parskip
 \bigskipamount=3\parskip
\jot=3pt
%
%
\def\ref#1{\par\noindent\hangindent2\parindent
 \hbox to 2\parindent{#1\hfil}\ignorespaces}
%
%
\font\typc=cmbx10 scaled \magstep1   
\font\typf=cmcsc10                   
\font\tenss=cmss10
\font\sevenss=cmss8 at 7pt
\font\fivess=cmss8 at 5pt
\newfam\ssfam %
\textfont\ssfam=\tenss
\scriptfont\ssfam=\sevenss
\scriptscriptfont\ssfam=\fivess
%
%
%
%
%
%
%
%
%
\catcode`\_=11
\def\suf_fix{}
\def\scaled_rm_box#1{%
 \relax
 \ifmmode
   \mathchoice
    {\hbox{\tenrm #1}}%
    {\hbox{\tenrm #1}}%
    {\hbox{\sevenrm #1}}%
    {\hbox{\fiverm #1}}%
 \else
  \hbox{\tenrm #1}%
 \fi}
\def\suf_fix_def#1#2{\expandafter\def\csname#1\suf_fix\endcsname{#2}}
\def\I_Buchstabe#1#2#3{%
 \suf_fix_def{#1}{\scaled_rm_box{I\hskip-0.#2#3em #1}}
}
\def\rule_Buchstabe#1#2#3#4{%
 \suf_fix_def{#1}{%
  \scaled_rm_box{%
   \hbox{%
    #1%
    \hskip-0.#2em%
    \lower-0.#3ex\hbox{\vrule height1.#4ex width0.07em }%
   }%
   \hskip0.50em%
  }%
 }%
}
\I_Buchstabe B22
\rule_Buchstabe C51{34}
\I_Buchstabe D22
\I_Buchstabe E22
\I_Buchstabe F22
\rule_Buchstabe G{525}{081}4
\I_Buchstabe H22
\I_Buchstabe I20
\I_Buchstabe K22
\I_Buchstabe L20
\I_Buchstabe M{20em }{I\hskip-0.35}
\I_Buchstabe N{20em }{I\hskip-0.35}
\rule_Buchstabe O{525}{095}{45}
\I_Buchstabe P20
\rule_Buchstabe Q{525}{097}{47}
\I_Buchstabe R21 
\rule_Buchstabe U{45}{02}{54}
\suf_fix_def{Z}{\scaled_rm_box{Z\hskip-0.38em Z}}
\catcode`\"=12
\newcount\math_char_code
\def\suf_fix_math_chars_def#1{%
 \ifcat#1A
  \expandafter\math_char_code\expandafter=\suf_fix_fam
  \multiply\math_char_code by 256
  \advance\math_char_code by `#1
  \expandafter\mathchardef\csname#1\suf_fix\endcsname=\math_char_code
  \let\next=\suf_fix_math_chars_def
 \else
  \let\next=\relax
 \fi
 \next}
%
%
%
%
\def\font_fam_suf_fix#1#2 #3 {%
 \def\suf_fix{#2}
 \def\suf_fix_fam{#1}
 \suf_fix_math_chars_def #3.
}
\font_fam_suf_fix
 0rm
 ABCDEFGHIJKLMNOPQRSTUVWXYZabcdefghijklmnopqrstuvwxyz
\font_fam_suf_fix
 2scr
 ABCDEFGHIJKLMNOPQRSTUVWXYZ
\font_fam_suf_fix
 \slfam sl
 ABCDEFGHIJKLMNOPQRSTUVWXYZabcdefghijklmnopqrstuvwxyz
\font_fam_suf_fix
 \bffam bf
 ABCDEFGHIJKLMNOPQRSTUVWXYZabcdefghijklmnopqrstuvwxyz
\font_fam_suf_fix
 \ttfam tt
 ABCDEFGHIJKLMNOPQRSTUVWXYZabcdefghijklmnopqrstuvwxyz
\font_fam_suf_fix
 \ssfam
 ss
 ABCDEFGHIJKLMNOPQRSTUVWXYZabcdefgijklmnopqrstuwxyz
\catcode`\_=8
\def\Adss{{\fam\ssfam A\mkern -9.5mu A}}%
\def\Cdss{{\fam\ssfam
    \mkern 4.2 mu \mathchoice%
    {\vrule height 6.5pt depth -.55pt width 1pt}%
    {\vrule height 6.5pt depth -.57pt width 1pt}%
    {\vrule height 4.55pt depth -.28pt width .8pt}%
    {\vrule height 3.25pt depth -.19pt width .6pt}%
    \mkern -6.3mu C}}%
\def\Ddss{{\fam\ssfam I\mkern -2.5mu D}}%
\def\Gdss{{\fam\ssfam
    \mkern 3.8 mu \mathchoice%
    {\vrule height 6.5pt depth -.62pt width 1pt}%
    {\vrule height 6.5pt depth -.65pt width 1pt}%
    {\vrule height 4.55pt depth -.44pt width .8pt}%
    {\vrule height 3.25pt depth -.30pt width .6pt}%
    \mkern -5.9mu G}}%
\def\Ndss{{\fam\ssfam I\mkern -2.5mu N}}%
\def\Qdss{{\fam\ssfam
    \mkern 3.8 mu \mathchoice%
    {\vrule height 6.5pt depth -.67pt width 1pt}%
    {\vrule height 6.5pt depth -.7pt width 1pt}%
    {\vrule height 4.55pt depth -.44pt width .7pt}%
    {\vrule height 3.25pt depth -.3pt width .5pt}%
    \mkern -5.9mu Q}}%
\def\Rdss{{\fam\ssfam I\mkern -2.5mu R}}%
\def\Zdss{{\fam\ssfam Z\mkern-8.1mu Z}}%
%
%
%
%
\font\teneuf=eufm10 
\font\seveneuf=eufm7
\font\fiveeuf=eufm5
\newfam\euffam \def\euf{\fam\euffam\teneuf} 
\textfont\euffam=\teneuf \scriptfont\euffam=\seveneuf
\scriptscriptfont\euffam=\fiveeuf

       \def\gfr{{\euf g}}

       \def\pfr{{\euf p}}

\input xypic.tex

\def\Hom{\mathop{\rm Hom}\nolimits}
\def\End{\mathop{\rm End}\nolimits}
\def\Dist{\mathop{\rm Dist}\nolimits}
\def\Lie{\mathop{\rm Lie}\nolimits}
\def\sup{\mathop{\rm sup}\nolimits}
\def\max{\mathop{\rm max}\nolimits}
\def\ind{\mathop{\rm ind}\nolimits}
\def\Ind{\mathop{\rm Ind}\nolimits}
\def\val{\mathop{\rm val}\nolimits}
\def\dlongrightarrow{\longrightarrow\hskip-8pt\rightarrow}

\font\typf=cmcsc10
\font\typff=cmcsc8

\parindent=0pt

\centerline{\typc Banach-Hecke algebras and $p$-adic Galois
representations}

\medskip

\centerline{\typf Peter Schneider and Jeremy Teitelbaum}

\bigskip

{\it Wir lassen vom Geheimnis uns erheben\hfill\break
 Der magischen Formelschrift, in deren Bann\hfill\break
 Das Uferlose, St\"urmende, das Leben\hfill\break
 Zu klaren Gleichnissen gerann.}\hfill\break
 \phantom\qquad\qquad{\typff Hermann Hesse}

\bigskip

\rightline{\it Dedicated to John Coates}

\bigskip

In this paper, we take some initial steps towards illuminating the
(hypothetical) $p$-adic local Langlands functoriality principle
relating Galois representations of a $p$-adic field $L$ and
admissible unitary Banach space representations of $G(L)$ when $G$
is a split reductive group over $L$.  The outline of our work is
derived from Breuil's remarkable insights into the nature of the
correspondence between $2$-dimensional crystalline Galois
representations of the Galois group of $\Qdss_p$ and Banach space
representations of $GL_{2}(\Qdss_p)$.

In the first part of the paper, we study the $p$-adic completion
$\Bscr(G,\rho)$ of the Hecke algebra $\Hscr(G,\rho)$ of
bi-equivariant compactly supported $\End(\rho)$-valued functions on
a totally disconnected, locally compact group $G$ derived from a
finite dimensional continuous representation $\rho$ of a compact
open subgroup $U$ of $G$. (These are the ``Banach-Hecke algebras''
of the title). After describing some general features of such
algebras we study in particular the case where $G$ is split
reductive and $U=U_0$ is a special maximal compact or $U=U_1$ is an
Iwahori subgroup of $G$ and $\rho$ is the restriction of a finite
dimensional algebraic representation of $G$ to $U_0$ or $U_1$.

In the smooth theory for trivial $\rho = 1_U$, by work of Bernstein,
the maximal commutative subalgebra of the Iwahori-Hecke algebra is
isomorphic to the group ring $K[\Lambda]$ where $\Lambda$ is the
cocharacter group of a maximal split torus $T$ of $G$, and the
spherical Hecke algebra is isomorphic by the Satake isomorphism to
the ring $K[\Lambda ]^{W}$ of Weyl group invariants.  At the same
time the algebra $K[\Lambda ]$ may be viewed as the ring of
algebraic functions on the dual maximal torus $T'$ in the dual group
$G'$. Together, these isomorphisms allow the identification of
characters of the spherical Hecke algebra with semisimple conjugacy
classes in $G'$. On the one hand, the Hecke character corresponds to
a certain parabolically induced smooth representation; on the other,
the conjugacy class in $G'$ determines the Frobenius in an
unramified Weil group representation of the field $L$.  This is the
unramified local Langlands correspondence (the Satake
parametrization) in the classical case.

With these principles in mind, we show that the completed maximal
commutative subalgebra of the Iwahori-Hecke algebra for $\rho$ is
isomorphic to the affinoid algebra of a certain explicitly given
rational subdomain $T'_{\rho}$ in the dual torus $T'$. The spectrum
of this algebra therefore corresponds to certain points of $T'$.  We
also show that the quotient of this subdomain by the Weyl group
action is isomorphic to the corresponding completion of the
spherical Hecke algebra; this algebra, for most groups $G$, turns
out to be a Tate algebra. These results may be viewed as giving a
$p$-adic completion of the Satake isomorphism, though our situation
is somewhat complicated by our reluctance to introduce a square root
of $q$ as is done routinely in the classical case. These
computations take up the first four sections of the paper.

In the second part of the paper, we let $G=GL_{d+1}(L)$.  We relate
the subdomain  of $T'$ determined by the completion $\Bscr(G,\rho)$
to isomorphism classes of a certain kind of crystalline Dieudonne
module.  This relationship follows Breuil's theory, which puts a
$2$-dimensional irreducible crystalline representation $V$ of ${\rm
Gal}(\overline{\Qdss_p}/\Qdss_p)$ with coefficients in a field $K$
into correspondence with a topologically irreducible admissible
unitary representation of $GL_2(\Qdss_p)$ in a $K$-Banach space.
Furthermore,  this Banach space representation is a completion of a
locally algebraic representation whose smooth factor comes from
$D_{cris}(V)$ viewed as a Weil group representation and whose
algebraic part is determined by the Hodge-Tate weights of $V$.

To state our relationship, let $V$ be a  $d+1$-dimensional
crystalline representation of ${\rm Gal}(\overline{L}/L)$ in a
$K$-vector space, where $K\supseteq L$ are finite extensions of
$\Qdss_p$. In this situation, $D_{cris}(V)$ has a $K$-vector space
structure. Suppose further that:

i. the eigenvalues of the Frobenius on $D_{cris}(V)$ lie in $K$;

ii. the (negatives of) the Hodge-Tate weights of $D_{cris}(V)$ are
multiplicity free and are separated from one another by at least
$[L:\Qdss_p]$;

iii. $V$ is {\it special}, meaning that the kernel of the natural
map
$$
\Cdss_p\otimes_{\Qdss_p} V\to \Cdss_{p}\otimes_{L} V
$$
is generated by its ${\rm Gal}(\overline{L}/L)$ invariants.

It follows from the Colmez-Fontaine theory that the category of such
special representations is equivalent to a category of
``$K$-isocrystals'', which are $K$-vector spaces with a $K$-linear
Frobenius and a filtration that is admissible in a sense very close
to the usual meaning.

Given such a representation, we extract from the associated
$K$-isocrystal its Frobenius, which we view as an element of the
dual group $G'(K)$ determined up to conjugacy.  The semi-simple part
$\zeta$ of this element determines a point of $T'(K)$ up to the Weyl
group action. From the Hodge-Tate weights, we extract a dominant
cocharacter of $G'$ and hence a highest $\xi$ determining an
algebraic representation $\rho = \rho_\xi$ for $G$. (In fact, the
highest weight is a modification of the Hodge-Tate weights, but we
avoid this complication in this introduction).   Put together, this
data yields a completion of the Iwahori-Hecke algebra, determined by
the highest weight, and a character of its maximal commutative
subalgebra, determined up to the Weyl group action. In other words,
we obtain a simple module $K_{\zeta}$ for the completed spherical
Hecke algebra $\Bscr(G,\rho_\xi|U_0)$.

Our main result is that  the existence of an admissible filtration
on $D_{cris}(V)$ translates into the condition that the point of
$T'$ determined by the Frobenius lives in the subdomain $T'_{\rho}$.
Conversely, we show how to reverse this procedure and, from a point
of $T'_{\rho}(K)$ (up to Weyl action), make an isocrystal that
admits an admissible filtration of Hodge-Tate type determined by
$\rho$. See Section 5 (esp.\ Proposition 5.2) for the details.

It is crucial to realize that the correspondence between points of
$T'_{\rho}$ and isocrystals outlined above does not determine a
specific filtration on the isocrystal.  Except when $d=1$ there are
infinitely many choices of filtration compatible with the given
data. Consequently the ``correspondence'' we describe is a very
coarse version of a $p$-adic local Langlands correspondence.

To better understand this coarseness on the
``representation-theoretic'' side, recall that to a Galois
representation $V$ of the type described above we associate a simple
module $K_{\zeta}$ for the completion $\Bscr(G,\rho|U_0)$ of the
spherical Hecke algebra. There is an easily described sup-norm on
the smooth compactly induced representation
$\ind_{U_0}^{G}(\rho|U_0)$; let  $B_{U_0}^G(\rho|U_0)$ be the
completion of this representation.  We show that the completed Hecke
algebra acts continuously on this space. By analogy with the
Borel-Matsumoto theory constructing parabolically induced
representations from compactly induced ones, and following also
Breuil's approach for $GL_{2}(\Qdss_p)$, it is natural to consider
the completed tensor product
$$
B_{\xi,\zeta} := K_{\zeta}\,
\widehat{\otimes}_{\Bscr(G,\rho_\xi|U_0)}\,
B_{U_0}^{G}(\rho_\xi|U_0)\ .
$$
A deep theorem of Breuil-Berger ([BB]) says that, in the
$GL_{2}(\Qdss_p)$-case, this representation in most cases is
nonzero, admissible, and irreducible, and under Breuil's
correspondence it is the Banach representation associated to $V$. In
our more general situation, we do not know even that $B_{\xi,\zeta}$
is nonzero. Accepting, for the moment, that it is nonzero, we do not
expect it to be admissible or irreducible, because it is associated
to the entire infinite family of representations having the same
Frobenius and Hodge-Tate weights as $V$ but different admissible
filtrations. We propose that $B_{\xi,\zeta}$ maps, with dense image,
to each of the Banach spaces coming from this family of Galois
representations. We discuss this further in Section 5.

In the last section of this paper (Section 6) we consider the shape
of a $p$-adic local Langlands functoriality for a general $L$-split
reductive group $G$ over $L$, with Langlands dual group $G'$ also
defined over $L$. Here we rely on ideas from the work of Kottwitz,
Rapoport-Zink, and Fontaine-Rapoport. Recall that a cocharacter
$\nu$ of the dual group $G'$ defined over $K$  allows one to put a
filtration $Fil^{\cdot}_{\rho'\circ\nu}\, E$ on every $K$-rational
representation space $(\rho',E)$ of $G'$. Using (a modified version
of) a notion of Rapoport-Zink, we say that a pair $(\nu,b)$
consisting of an element $b$ of $G'(K)$ and a $K$-rational
cocharacter $\nu$ of $G'$  is an ``admissible pair'' if, for any
$K$-rational representation $(\rho',E)$ of $G'$, the $K$-isocrystal
$(E,\rho'(b),Fil^\cdot_{\rho'\circ\nu}\, E)$ is admissible. Such an
admissible pair defines a faithful tensor functor from the neutral
Tannakian category of $K$-rational representations of $G'$ to that
of admissible filtered $K$-isocrystals.  Composing this with the
Fontaine functor one obtains a tensor functor to the category of
``special'' ${\rm Gal}(\overline{L}/L)$ representations of the type
described earlier. The Tannakian formalism therefore constructs from
an admissible pair an isomorphism class of representations of the
Galois group of $L$ in $G'(\overline{K})$.

Now suppose given an irreducible algebraic representation $\rho$ of
$G$. Its highest weight  may be viewed as a cocharacter of $G'$.
Under a certain technical condition, we prove in this section that
there is an admissible pair $(\nu,b)$ where $\nu$ is conjugate by
$G'(K)$ to a (modification of) the highest weight, and $b$ is an
element of $G'(K)$, if and only if the semisimple part of $b$ is
conjugate to an element of the affinoid domain $T'_{\rho}(K)$ (See
Proposition 6.1). Thus in some sense this domain is functorial in
the group $G'$.

Our work in this section relies on a technical hypothesis on $G$.
Suppose that $\eta$ is half the sum of the positive roots of $G$. We
need $[L:\Qdss_p]\eta$ to be an integral weight of $G$.  This
happens, for example, if $L$ has even degree over $\Qdss_p$, and in
general for many groups, but not, for example, when
$G=PGL_{2}(\Qdss_p)$. This complication has its origin in the
normalization of the Langlands correspondence. Because of the square
root of $q$ issue the $p$-adic case seems to force the use of the
``Hecke'' or the ``Tate'' correspondence rather than the traditional
unitary correspondence; but even for smooth representations this is
not functorial (cf.\ [Del] (3.2.4-6)). It turns out that without the
above integrality hypothesis one even has to introduce a square root
of a specific continuous Galois character (for $L = \Qdss_p$ it is
the cyclotomic character). This leads to isocrystals with a
filtration indexed by half-integers. Although it seems possible to
relate these to Galois representations this has not been done yet in
the literature. We hope to come back to this in the future.

The authors thank Matthew Emerton for pointing out that the
conditions which define our affinoid domains $T'_\rho$ are
compatible with the structure of his Jacquet functor on locally
algebraic representations ([Em1] Prop.\ 3.4.9 and Lemma 4.4.2, [Em2]
Lemma 1.6). We thank Laurent Berger, Christophe Breuil, and
especially Jean-Marc Fontaine for their very helpful conversations
about these results. We also want to stress that our computations in
Section 4 rely in an essential way on the results of Marie-France
Vigneras in [Vig]. The first author gratefully acknowledges support
from UIC and CMI. During the final stages of this paper he was
employed by the Clay Mathematics Institute as a Research Scholar.
The second author was supported by National Science Foundation Grant
DMS-0245410.

We dedicate this paper to John Coates on the occasion of his
sixtieth birthday. His constant support and unrelenting enthusiasm
was and is an essential source of energy and inspiration for us over
all these years.

Throughout this paper $K$ is a fixed complete extension field of
$\Qdss_p$ with absolute value $|\ |$.

\medskip

{\bf 1. Banach-Hecke algebras}

\medskip

In this section $G$ denotes a totally disconnected and locally
compact group, and $U \subseteq G$ is a fixed compact open
subgroup. We let $(\rho,W)$ be a continuous representation of $U$
on a finite dimensional $K$-vector space $W$, and we fix a
$U$-invariant norm $\|\ \|$ on $W$.

The Hecke algebra $\Hscr(G,\rho)$ is the $K$-vector space of all
compactly supported functions $\psi : G \longrightarrow \End_K(W)$
satisfying
$$
\psi(u_1gu_2) = \rho(u_1)\circ \psi(g) \circ \rho(u_2)
\qquad\hbox{for any}\ u_1,u_2 \in U\ \hbox{and}\ g \in G\ .
$$
It is a unital associative $K$-algebra via the convolution
$$
\psi_1 \ast \psi_2(h) := \sum_{g \in G/U} \psi_1(g) \circ
\psi_2(g^{-1}h)\ .
$$
Its unit element is the function
$$
\psi_e(h) := \left\{\matrix{\rho(h) & \hbox{if}\ h \in U,\hfill\cr
0 \hfill & \hbox{otherwise}.}\right.
$$
We note that any function $\psi$ in $\Hscr(G,\rho)$ necessarily is
continuous. We now introduce the norm
$$
\|\psi\| := \sup_{g \in G} \|\psi(g)\|
$$
on $\Hscr(G,\rho)$ where on the right hand side $\|\ \|$ refers to
the operator norm on $\End_K(W)$ with respect to the original norm
$\|\ \|$ on $W$. This norm on $\Hscr(G,\rho)$ evidently is
submultiplicative. By completion we therefore obtain a unital
$K$-Banach algebra $\Bscr(G,\rho)$, called in the following the
Banach-Hecke algebra, with submultiplicative norm. As a Banach
space $\Bscr(G,\rho)$ is the space of all continuous functions
$\psi : G \longrightarrow \End_K(W)$ vanishing at infinity and
satisfying
$$
\psi(u_1gu_2) = \rho(u_1)\circ \psi(g) \circ \rho(u_2)
\qquad\hbox{for any}\ u_1,u_2 \in U\ \hbox{and}\ g \in G\ .
$$
In the special case where $\rho = 1_U$ is the trivial
representation $\Hscr(G,1_U)$, resp. $\Bscr(G,1_U)$, is the vector
space of all $K$-valued finitely supported functions, resp.
functions vanishing at infinity, on the double coset space
$U\backslash G/U$.

A more intrinsic interpretation of these algebras can be given by
introducing the compactly induced $G$-representation
$\ind_U^G(\rho)$. This is the $K$-vector space of all compactly
supported functions $f : G \longrightarrow W$ satisfying
$$
f(gu) = \rho(u^{-1})(f(g))\qquad\hbox{for any}\ u \in U\
\hbox{and}\ g \in G
$$
with $G$ acting by left translations. Again we note that any
function $f$ in $\ind_U^G(\rho)$ is continuous. We equip
$\ind_U^G(\rho)$ with the $G$-invariant norm
$$
\|f\| := \sup_{g \in G} \|f(g)\|
$$
and let $B_U^G(\rho)$ denote the corresponding completion. The
$G$-action extends isometrically to the $K$-Banach space
$B_U^G(\rho)$, which consists of all continuous functions $f : G
\longrightarrow W$ vanishing at infinity and satisfying
$$
f(gu) = \rho(u^{-1})(f(g))\qquad\hbox{for any}\ u \in U\
\hbox{and}\ g \in G\ .
$$

\smallskip

{\bf Lemma 1.1:} {\it The $G$-action on $B_U^G(\rho)$ is
continuous.}

Proof: Since $G$ acts isometrically it remains to show that the
orbit maps
$$
\matrix{ c_f : G & \longrightarrow & B_U^G(\rho)\cr
 \hfill g & \longmapsto & gf\ , \hfill}
$$
for any $f \in B_U^G(\rho)$, are continuous. In case $f \in
\ind_U^G(\rho)$ the map $c_f$ even is locally constant. In general
we write $f = \mathop{\rm lim}\limits_{i \rightarrow \infty}$ as
the limit of a sequence $(f_i)_{i \in \Ndss}$ in $\ind_U^G(\rho)$.
Because of
$$
\|(c_f - c_{f_i})(g)\| = \|g(f - f_i)\| = \|f - f_i\|
$$
the map $c_f$ is the uniform limit of the locally constant maps
$c_{f_i}$ and hence is continuous.

\smallskip

One easily checks that the pairing
$$
\matrix{
 \Hscr(G,\rho) \times \ind_U^G(\rho) & \longrightarrow &
\ind_U^G(\rho) \hfill\cr \hfill (\psi,f) & \longmapsto & (\psi\ast
f)(g) := \sum_{h \in G/U} \psi(g^{-1}h)(f(h)) } \leqno{(1)}
$$
makes $\ind_U^G(\rho)$ into a unital left $\Hscr(G,\rho)$-module
and that this module structure commutes with the $G$-action.

\medskip

{\bf Lemma 1.2:} {\it The map
$$
\matrix{
 \Hscr(G,\rho) & \mathop{\longrightarrow}\limits^{\cong} &
 \End_G(\ind_U^G(\rho)) \hfill\cr
 \hfill \psi & \longmapsto & A_{\psi}(f) := \psi \ast f \hfill }
$$
is an isomorphism of $K$-algebras.}

Proof: For a smooth representation $\rho$ this can be found in
[Kut]. Our more general case follows by the same argument. But
since we will need the notations anyway we recall the proof. The
map in question certainly is a homomorphism of $K$-algebras. We
now introduce, for any $w \in W$, the function
$$
f_w(g) := \left\{\matrix{\rho(g^{-1})(w) & \hbox{if}\ g \in
U,\hfill\cr 0 \hfill & \hbox{otherwise}}\right.
$$
in $\ind_U^G(\rho)$. We have
$$
A_{\psi}(f_w)(g) = (\psi\ast f_w)(g) = \psi(g^{-1})(w)\quad
\hbox{for any}\ \psi \in \Hscr(G,\rho)\ .\leqno{(2)}
$$
This shows that the map in question is injective. To see its
surjectivity we fix an operator $A_0 \in \End_G(\ind_U^G(\rho))$
and consider the function
$$
\matrix{
 \psi_0 : G & \longrightarrow & \End_K(W) \hfill\cr
 \hfill g & \longmapsto & [w \mapsto A_0(f_w)(g^{-1})]\ .}
$$
It clearly has compact support. Furthermore, for $u_1,u_2 \in U$,
we compute
$$
\eqalign{
 \psi_0(u_1gu_2)(w) & = A_0(f_w)(u_2^{-1}g^{-1}u_1^{-1})
 = \rho(u_1)[A_0(f_w)(u_2^{-1}g^{-1})] \cr
 & = \rho(u_1)[(u_2(A_0(f_w)))(g^{-1})] =
 \rho(u_1)[A_0(u_2(f_w))(g^{-1})] \cr
 & = \rho(u_1)[A_0(f_{\rho(u_2)(w)})(g^{-1})] =
 \rho(u_1)[\psi_0(\rho(u_2)(w))] \cr
 & = [\rho(u_1) \circ \psi_0 \circ \rho(u_2)](w)\ . }
$$
Hence $\psi_0 \in \Hscr(G,\rho)$. Moreover, for any $f \in
\ind_U^G(\rho)$ we have
$$
f = \sum_{h \in G/U} h(f_{f(h)})
$$
and therefore
$$
\eqalign{
 A_{\psi_0}(f)(g) & = (\psi_0\ast f)(g) = \sum_{h \in G/U}
 \psi_0(g^{-1}h)(f(h)) \cr
 & = \sum_{h \in G/U} A_0(f_{f(h)})(h^{-1}g) = A_0(\sum_{h \in G/U}
 h(f_{f(h)}))(g) \cr
 & = A_0(f)(g)\ .}
$$
Hence $A_{\psi_0} = A_0$.

\medskip

We evidently have $\|\psi \ast f\| \leq \|\psi\| \cdot \|f\|$. By
continuity we therefore obtain a continuous left action of the
Banach algebra $\Bscr(G,\rho)$ on the Banach space $B_U^G(\rho)$
which is submultiplicative in the corresponding norms and which
commutes with the $G$-action. This action is described by the same
formula (1), and we therefore continue to denote it by $\ast$.

\medskip

{\bf Lemma 1.3:} {\it The map
$$
\matrix{
 \Bscr(G,\rho) & \mathop{\longrightarrow}\limits^{\cong} &
 \End_G^{cont}(B_U^G(\rho)) \hfill\cr
 \hfill \psi & \longmapsto & A_{\psi}(f) := \psi \ast f \hfill }
$$
is an isomorphism of $K$-algebras and is an isometry with respect
to the operator norm on the right hand side.}

Proof: (The superscript ``cont'' refers to the continuous
endomorphisms.) By the previous discussion the map $\psi
\longmapsto A_{\psi}$ is well defined, is a homomorphism of
$K$-algebras, and is norm decreasing. Using the notations from the
proof of Lemma 1.2 the formula (2), by continuity, holds for any
$\psi \in \Bscr(G,\rho)$. Using that $\|f_w\| = \|w\|$ we now
compute
$$
\eqalign{
 \|A_{\psi}\| & \geq \mathop{\rm sup}\limits_{w \neq 0}
 {\|\psi\ast f_w\| \over \|f_w\|} =
 \mathop{\rm sup}\limits_{w \neq 0} \mathop{\rm sup}\limits_g
 {\|\psi(g^{-1})(w)\|\over \|w\|} = \mathop{\rm sup}\limits_g
 \|\psi(g^{-1})\| \cr
 & = \|\psi\| \geq \|A_\psi\| \ .}
$$
It follows that the map in the assertion is an isometry and in
particular is injective. To see its surjectivity we fix an $A_0
\in \End_G^{cont}(B_U^G(\rho))$ and define
$$
\matrix{
 \psi_0 : G & \longrightarrow & \End_K(W) \hfill\cr
 \hfill g & \longmapsto & [w \mapsto A_0(f_w)(g^{-1})] \ .}
$$
Since each $A_0(f_w)$ is continuous and vanishing at infinity on
$G$ it follows that $\psi_0$ is continuous and vanishing at
infinity. By exactly the same computations as in the proof of
Lemma 1.2 one then shows that in fact $\psi_0 \in \Bscr(G,\rho)$
and that $A_{\psi_0} = A_0$.

\medskip

We end this section by considering the special case where
$(\rho,E)$ is the restriction to $U$ of a continuous
representation $\rho$ of $G$ on a finite dimensional $K$-vector
space $E$. It is easy to check that then the map
$$
\matrix{
 \iota_\rho : \Hscr(G,1_U) & \longrightarrow & \Hscr(G,\rho) \cr
 \hfill \psi & \longmapsto & \psi\cdot\rho \hfill}
$$
is an injective homomorphism of $K$-algebras. There are
interesting situations where this map in fact is an isomorphism.
We let $L$ be a finite extension of $\Qdss_p$ contained in $K$,
and we assume that $G$ as well as $(\rho,E)$ are locally
$L$-analytic.

\medskip

{\bf Lemma 1.4:} {\it Suppose that, for the derived action of the
Lie algebra $\gfr$ of $G$, the $K \otimes_L \gfr$-module $E$ is
absolutely irreducible; then the homomorphism $\iota_\rho$ is
bijective.}

Proof: Using Lemma 1.2 and Frobenius reciprocity we have
$$
\eqalign{
 \Hscr(G,\rho) & = \End_G(\ind^G_U(\rho)) =
 \Hom_U(E,\ind^G_U(\rho)) \cr
 & = \Hom_U(E,\ind^G_U(1) \otimes_K E) \cr
 & = [\ind^G_U(1) \otimes_K E^\ast \otimes_K E]^U }
$$
where the last term denotes the $U$-fixed vectors in the tensor
product with respect to the diagonal action. This diagonal action
makes the tensor product equipped with the finest locally convex
topology into a locally analytic $G$-representation. Its $U$-fixed
vectors certainly are contained in the vectors annihilated by the
derived action of $\gfr$. Since $G$ acts smoothly on $\ind^G_U(1)$
we have
$$
\eqalign{
 (\ind^G_U(1) \otimes_K E^\ast \otimes_K E)^{\gfr = 0} & =
 \ind^G_U(1) \otimes_K (E^\ast \otimes_K E)^{\gfr = 0} \cr
 & = \ind^G_U(1) \otimes_K \End_{K \otimes_L \gfr}(E) \ .}
$$
Our assumption on absolute irreducibility implies that $\End_{K
\otimes_L \gfr}(E) = K$. We therefore see that
$$
\Hscr(G,\rho) = [\ind^G_U(1) \otimes_K E^\ast \otimes_K E]^U =
\ind^G_U(1) = \Hscr(G,1_U) \ .
$$

\medskip

{\bf 2. Weights and affinoid algebras}

\medskip

For the rest of this paper $L/\Qdss_p$ is a finite extension
contained in $K$, and $G$ denotes the group of $L$-valued points
of an $L$-split connected reductive group over $L$. Let $|\ |_L$
be the normalized absolute value of $L$, $\val_L : K^\times
\longrightarrow \Rdss$ the unique additive valuation such that
$\val_L(L^\times) = \Zdss$, and $q$ the number of elements in the
residue class field of $L$. We fix a maximal $L$-split torus $T$
in $G$ and a Borel subgroup $P = TN$ of $G$ with Levi component
$T$ and unipotent radical $N$. The Weyl group of $G$ is the
quotient $W = N(T)/T$ of the normalizer $N(T)$ of $T$ in $G$ by
$T$. We also fix a maximal compact subgroup $U_0 \subseteq G$
which is special with respect to $T$ (i.e., is the stabilizer of a
special vertex $x_0$ in the apartment corresponding to $T$, cf.
[Car]\S3.5). We put $T_0 := U_0 \cap T$ and $N_0 := U_0 \cap N$.
The quotient $\Lambda := T/T_0$ is a free abelian group of rank
equal to the dimension of $T$ and can naturally be identified with
the cocharacter group of $T$. Let $\lambda : T \longrightarrow
\Lambda$ denote the projection map. The conjugation action of
$N(T)$ on $T$ induces $W$-actions on $T$ and $\Lambda$ which we
denote by $t \longmapsto {^w}t$ and $\lambda \longmapsto
{^w}\lambda$, respectively. We also need the $L$-torus $T'$ dual
to $T$. Its $K$-valued points are given by
$$
T'(K) := \Hom(\Lambda,K^\times) \ .
$$
The group ring $K[\Lambda]$ of $\Lambda$ over $K$ naturally
identifies with the ring of algebraic functions on the torus $T'$.
We introduce the ``valuation map''
$$
\xymatrix{
  val :\ T'(K) = \Hom(\Lambda,K^\times) \ar[rr]^{\quad\val_L\circ} & &
  \Hom(\Lambda,\Rdss) =: V_\Rdss  \ . }
$$
If $X^\ast(T)$ denotes the algebraic character group of the torus
$T$ then we have the embedding
$$
\matrix{
 X^\ast(T) & \longrightarrow & \Hom(\Lambda,\Rdss) \cr
 \hfill \chi & \longmapsto & \val_L\circ\chi \hfill }
$$
which induces an isomorphism
$$
X^\ast(T) \otimes \Rdss \mathop{\rightarrow}\limits^{\cong}
V_\Rdss \ .
$$
We therefore may view $V_\Rdss$ as the real vector space
underlying the root datum of $G$ with respect to $T$. Evidently
any $\lambda \in \Lambda$ defines a linear form in the dual vector
space $V_\Rdss^\ast$ also denoted by $\lambda$. Let $\Phi$ denote
the set of roots of $T$ in $G$ and let $\Phi^+ \subseteq \Phi$ be
the subset of those roots which are positive with respect to $P$.
As usual, $\check{\alpha} \in \Lambda$ denotes the coroot
corresponding to the root $\alpha \in \Phi$. The subset
$\Lambda^{--} \subseteq \Lambda$ of antidominant cocharacters is
defined to be the image $\Lambda^{--} := \lambda(T^{--})$ of
$$
T^{--} := \{t \in T : |\alpha(t)|_L \geq 1\ \hbox{for any}\ \alpha
\in \Phi^+\} \ .
$$
Hence
$$
\Lambda^{--} = \{ \lambda \in \Lambda : \val_L \circ
\alpha(\lambda) \leq 0\ \hbox{for any}\ \alpha \in \Phi^+ \} \ .
$$
We finally recall that $\Lambda^{--}$ carries the partial order
$\leq$ given by
$$
\mu \leq \lambda \quad\hbox{if}\quad \lambda - \mu \in
\sum_{\alpha \in \Phi^+} \Rdss_{\geq 0} \cdot (- \check{\alpha})
\subseteq \Lambda \otimes \Rdss \ .
$$

In this section we will investigate certain Banach algebra
completions of the group ring $K[\Lambda]$ together with certain
twisted $W$-actions on them. We will proceed in an axiomatic way
and will give ourselves a cocycle on $W$ with values in $T'(K)$,
i.e., a map
$$
\gamma : W \times \Lambda \longrightarrow K^\times
$$
such that
$$
\gamma(w,\lambda\mu) = \gamma(w,\lambda)\gamma(w,\mu)\qquad
\hbox{for any}\ w \in W\ \hbox{and}\ \lambda,\mu \in \Lambda
\leqno{(a)}
$$
and
$$
\gamma(vw,\lambda) = \gamma(v,{^w}\lambda)\gamma(w,\lambda)\qquad
\hbox{for any}\ v,w \in W\ \hbox{and}\ \lambda \in \Lambda \
.\leqno{(b)}
$$
Moreover we impose the positivity condition
$$
|\gamma(w,\lambda)| \leq 1\qquad \hbox{for any}\ w \in W\
\hbox{and}\ \lambda \in \Lambda^{--} \leqno{(c)}
$$
as well as the partial triviality condition
$$
\gamma(w,\lambda) = 1\qquad \hbox{for any}\ w \in W\ \hbox{and}\
\lambda \in \Lambda\ \hbox{such that}\ {^w}\lambda = \lambda \
.\leqno{(d)}
$$
The twisted action of $W$ on $K[\Lambda]$ then is defined by
$$
\matrix{
 \hfill W \times K[\Lambda] & \longrightarrow & K[\Lambda] \hfill\cr
 \hfill (w, \sum_\lambda c_\lambda \lambda) & \longmapsto &
 w\cdot(\sum_\lambda c_\lambda \lambda) := \sum_\lambda
 \gamma(w,\lambda)c_\lambda {^w\lambda}\ .}
$$
By $(a)$, each $w \in W$ acts as an algebra automorphism, and the
cocycle condition $(b)$ guarantees the associativity of this
$W$-action. The invariants with respect to this action will be
denoted by $K[\Lambda]^{W,\gamma}$. Since $\Lambda^{--}$ is a
fundamental domain for the $W$-action on $\Lambda$ it follows that
$K[\Lambda]^{W,\gamma}$ has the $K$-basis
$\{\sigma_\lambda\}_{\lambda \in \Lambda^{--}}$ defined by
$$
\sigma_\lambda := \sum_{w \in W/W(\lambda)} w\cdot\lambda =
\sum_{w \in W/W(\lambda)} \gamma(w,\lambda) {^w\lambda}
$$
where $W(\lambda) \subseteq W$ denotes the stabilizer of $\lambda$
and where the sums are well defined because of $(d)$. Next, again
using $(d)$, we define the map
$$
\matrix{
 \gamma^{dom} : \Lambda & \longrightarrow & K^\times \hfill\cr
 \hfill \lambda & \longmapsto & \gamma(w,\lambda)\quad \hbox{if}\
 \ {^w\lambda} \in \Lambda^{--} \ , }
$$
and we equip $K[\Lambda]$ with the norm
$$
\|\sum_\lambda c_\lambda \lambda\|_\gamma :=
\mathop{\sup}\limits_{\lambda \in \Lambda}
|\gamma^{dom}(\lambda)c_\lambda|\ .
$$
The cocycle condition $(b)$ implies the identity
$$
\gamma^{dom}({^w\lambda})\gamma(w,\lambda) =
\gamma^{dom}(\lambda)\leqno{(1)}
$$
from which one deduces that the twisted $W$-action on $K[\Lambda]$
is isometric in the norm $\|\ \|_\gamma$ and hence extends by
continuity to a $W$-action on the completion
$K\langle\Lambda;\gamma\rangle$ of $K[\Lambda]$ with respect to
$\|\ \|_\gamma$. Again we denote the corresponding $W$-invariants
by $K\langle\Lambda;\gamma\rangle^{W,\gamma}$. One easily checks
that $\{\sigma_\lambda\}_{\lambda \in \Lambda^{--}}$ is an
orthonormal basis of the Banach space
$(K\langle\Lambda;\gamma\rangle^{W,\gamma},\|\ \|_\gamma)$.

\medskip

{\bf Lemma 2.1:} {\it i. $|\gamma^{dom}(\lambda)| \geq 1$ for any
$\lambda \in \Lambda$;

ii. $|\gamma^{dom}(\lambda\mu)| \leq
|\gamma^{dom}(\lambda)||\gamma^{dom}(\mu)|$ for any $\lambda, \mu
\in \Lambda$. }

Proof: i. If $^w\lambda \in \Lambda^{--}$ then
$\gamma^{dom}(\lambda) = \gamma(w,\lambda) =
\gamma(w^{-1},{^w\lambda})^{-1}$. The claim therefore is a
consequence of the positivity condition $(c)$. ii. If
$^w(\lambda\mu) \in \Lambda^{--}$ then, using (1), we have
$$
\gamma^{dom}(\lambda\mu) =
\gamma^{dom}({^w\lambda})^{-1}\gamma^{dom}({^w\mu})^{-1}\gamma^{dom}(\lambda)
\gamma^{dom}(\mu) \ .
$$
Hence the claim follows from the first assertion.

\medskip

It is immediate from Lemma 2.1.i that the norm $\|\ \|_\gamma$ is
submultiplicative. Hence $K\langle\Lambda;\gamma\rangle$ is a
$K$-Banach algebra containing $K[\Lambda]$ as a dense subalgebra.
Moreover, since the twisted $W$-action on
$K\langle\Lambda;\gamma\rangle$ is by algebra automorphisms,
$K\langle\Lambda;\gamma\rangle^{W,\gamma}$ is a Banach subalgebra
of $K\langle\Lambda;\gamma\rangle$.

In order to compute the Banach algebra
$K\langle\Lambda;\gamma\rangle$ we introduce the subset
$$
T'_\gamma(K) := \{ \zeta \in T'(K) : |\zeta(\lambda)| \leq
|\gamma^{dom}(\lambda)|\ \hbox{for any}\ \lambda \in \Lambda \}
$$
of $T'(K)$. We obviously have
$$
T'_\gamma(K) = val^{-1}(V_\Rdss^\gamma)
$$
with
$$
V_\Rdss^\gamma := \{ z \in V_\Rdss : \lambda(z) \geq
\val_L(\gamma^{dom}(\lambda))\ \hbox{for any}\ \lambda \in
\Lambda\} \ .
$$
By $(a)$, our cocycle $\gamma$ defines the finitely many points
$$
z_w := - val(\gamma(w^{-1},.)) \qquad\hbox{for}\ w \in W
$$
in $V_\Rdss$. The cocycle condition $(b)$ implies that
$$
z_{vw} = {^v}z_w + z_v \qquad\hbox{for any}\ v,w \in W \leqno{(2)}
$$
and the positivity condition $(c)$ that
$$
\lambda(z_w) \leq 0 \qquad\hbox{for any}\ w \in W\ \hbox{and}\
\lambda \in \Lambda^{--} \ . \leqno{(3)}
$$

\medskip

{\bf Remark 2.2:} $\{z \in V_\Rdss : \lambda(z) \leq 0\ \hbox{for
any}\ \lambda \in \Lambda^{--}\} = \sum_{\alpha \in \Phi^+}
\Rdss_{\geq 0} \cdot \val_L\circ\alpha$.

Proof: This reduces to the claim that the (closed) convex hull of
$\Lambda^{--}$ in $V_\Rdss^\ast$ is equal to the antidominant cone
$$
(V_\Rdss^\ast)^{--} = \{ z^\ast \in V_\Rdss^\ast : z^\ast(z) \leq
0\ \hbox{for any}\ z \in \sum_{\alpha \in \Phi^+} \Rdss_{\geq 0}
\cdot \val_L\circ \alpha\} \ .
$$
Let $Z \subseteq G$ denote the connected component of the center
of $G$. Then $G/Z$ is semisimple and the sequence
$$
0 \longrightarrow Z/Z_0 \longrightarrow T/T_0 \longrightarrow
(T/Z)/(T/Z)_0 \longrightarrow 0
$$
is exact. Hence the fundamental antidominant coweights for the
semisimple group $G/Z$ can be lifted to elements
$\omega_1,\ldots,\omega_d \in V_\Rdss^\ast$ in such a way that,
for some $m \in \Ndss$, we have $m\omega_1,\ldots, m\omega_d \in
\Lambda^{--}$. It follows that
$$
(V_\Rdss^\ast)^{--} = (Z/Z_0) \otimes \Rdss + \sum_{i=1}^d
\Rdss_{\geq 0} \cdot \omega_i
$$
and
$$
\Lambda^{--} \supseteq Z/Z_0 + m \cdot \sum_{i=1}^d \Zdss_{\geq 0}
\cdot \omega_i \ .
$$

\medskip

We therefore obtain from $(3)$ that
$$
z_w \in \sum_{\alpha \in \Phi^+} \Rdss_{\geq 0} \cdot
\val_L\circ\alpha \qquad\hbox{for any}\ w \in W \ . \leqno{(4)}
$$
In terms of these points $z_w$ we have
$$
\eqalign{
 V_\Rdss^\gamma &
 = \{z \in V_\Rdss : \lambda(z) \geq \lambda(-z_{w^{-1}})\
 \hbox{for any}\ \lambda \in \Lambda, \ w \in W\
 \hbox{such that}\ {^w}\lambda \in \Lambda^{--}\} \cr
 &
 = \{z \in V_\Rdss : {^{w^{-1}}}\lambda(z) \geq {^{w^{-1}}}\lambda(-z_{w^{-1}})\
 \hbox{for any}\ w \in W\
 \hbox{and}\ \lambda \in \Lambda^{--}\} \cr
 &
 = \{z \in V_\Rdss : \lambda({^w}z) \geq \lambda(z_w)\
 \hbox{for any}\ w \in W\
 \hbox{and}\ \lambda \in \Lambda^{--}\}
 }
$$
where the last identity uses $(2)$. Obviously $V_\Rdss^\gamma$ is
a convex subset of $V_\Rdss$. Using the partial order $\leq$ on
$V_\Rdss$ defined by $\Phi^+$ (cf.\ [B-GAL] Chap.\ VI \S1.6) we
obtain from Remark 2.2 that
$$
V_\Rdss^\gamma = \{ z \in V_\Rdss : {^w}z \leq z_w\ \hbox{for
any}\ w \in W\} \ .
$$

\medskip

{\bf Lemma 2.3:} {\it $V_\Rdss^\gamma$ is the convex hull in
$V_\Rdss$ of the finitely many points $-z_w$ for $w \in W$.}

Proof: From $(2)$ and $(4)$ we deduce that
$$
^{w}z_v + z_w = z_{wv} \geq 0 \quad\hbox{and hence}\quad
^{w}(-z_v) \leq z_w
$$
for any $v,w \in W$. It follows that all $-z_v$ and therefore
their convex hull is contained in $V_\Rdss^\gamma$. For the
reverse inclusion suppose that there is a point $z \in
V_\Rdss^\gamma$ which does not lie in the convex hull of the
$-z_w$. We then find a linear form $\ell \in V_\Rdss^\ast$ such
that $\ell(z) < \ell(-z_w)$ for any $w \in W$. Choose $v \in W$
such that $\ell_0 := {^v}\ell$ is antidominant. It follows that
${^{v^{-1}}}\ell_0(z) < {^{v^{-1}}}\ell_0(-z_w)$ and hence, using
$(2)$, that
$$
\ell_0({^v}z) < \ell_0(-{^v}z_w) = \ell_0(z_v) - \ell_0(z_{vw})
$$
for any $w \in W$. For $w := v^{-1}$ we in particular obtain
$$
\ell_0({^v}z) < \ell_0(z_v) \ .
$$
On the other hand, since $z \in V_\Rdss^\gamma$, we have
$$
\lambda({^v}z) \geq \lambda(z_v)
$$
for any $\lambda \in \Lambda^{--}$ and hence for any $\lambda$ in
the convex hull of $\Lambda^{--}$. But as we have seen in the
proof of Remark 2.2 the antidominant $\ell_0$ belongs to this
convex hull which leads to a contradiction.

\medskip

{\bf Proposition 2.4:} {\it i. $T'_\gamma(K)$ is the set of
$K$-valued points of an open $K$-affinoid subdomain $T'_\gamma$ in
the torus $T'$;

ii. the Banach algebra $K\langle \Lambda;\gamma \rangle$ is
naturally isomorphic to the ring of analytic functions on the
affinoid domain $T'_\gamma$;

iii. the affinoid domain $T'_\gamma$ is the preimage, under the
map ``$val$'', of the convex hull of the finitely many points $-
z_w \in V_\Rdss$ for $w \in W$;

iv. $K\langle \Lambda;\gamma \rangle^{W,\gamma}$ is an affinoid
$K$-algebra. }

Proof: It follows from Gordan's lemma ([KKMS] p.\ 7) that the
monoid $\Lambda^{--}$ is finitely generated.  Choose a finite set
of generators $F^{--}$, and let
$$
F:= \{^{w}\lambda : \lambda\in F^{--}\} \ .
$$
Using the fact that, by construction, the function $\gamma^{dom}$
is multiplicative within Weyl chambers we see that the infinitely
many inequalities implicit in the definition of $T'_\gamma(K)$ can
in fact be replaced by finitely many:
$$
T'_\gamma(K)=\{\zeta\in T'(K) : |\zeta(\lambda)|\le
|\gamma^{dom}(\lambda)| \hbox{\rm\ for any $\lambda\in F$}\}.
$$
We therefore define $T'_\gamma$ to be the rational subset in $T'$
given by the finitely many inequalities
$|\gamma^{dom}(\lambda)^{-1}\lambda(\zeta)| \leq 1$ for $\lambda
\in F$ and obtain point i. of our assertion. Now choose
indeterminates $T_{\lambda}$ for $\lambda\in F$ and consider the
commutative diagram of algebra homomorphisms
$$
\xymatrix{
  o_K[T_\lambda : \lambda \in F] \ar[d]_{\subseteq} \ar[r]
       & K[\Lambda]^0 \ar[d]^{\subseteq} \\
  K[T_\lambda : \lambda \in F] \ar[d]_{\subseteq} \ar[r]
       & K[\Lambda] \ar[d]^{\subseteq} \\
  K\langle T_{\lambda}:\lambda\in F\rangle \ar[r]
       & K\langle\Lambda;\gamma\rangle  }
$$
where the horizontal arrows send  $T_{\lambda}$ to
$\gamma^{dom}(\lambda)^{-1}\lambda$, where $o_K$ is the ring of
integers of $K$, and where $K[\Lambda]^0$ denotes the unit ball
with respect to $\|\ \|_\gamma$ in $K[\Lambda]$. Again, the
multiplicativity of $\gamma^{dom}$ within Weyl chambers shows that
all three horizontal maps are surjective.  The lower arrow gives a
presentation of $K\langle\Lambda;\gamma\rangle$ as an affinoid
algebra. The middle arrow realizes the dual torus $T'$ as a closed
algebraic subvariety
$$
\matrix{
 \iota : T' & \longrightarrow & \Adss^f \hfill\cr
 \hfill \zeta & \longmapsto & (\zeta(\lambda))_{\lambda \in F} }
$$
in the affine space $\Adss^f$ where $f$ denotes the cardinality of
the set $F$. The surjectivity of the upper arrow shows that the
norm $\|\ \|_\gamma$ on $K[\Lambda]$ is the quotient norm of the
usual Gauss norm on the polynomial ring $K[T_\lambda : \lambda \in
F]$. Hence the kernel of the lower arrow is the norm completion of
the kernel $I$ of the middle arrow. Since any ideal in the Tate
algebra $K\langle T_{\lambda}:\lambda\in F\rangle$ is closed we
obtain
$$
K\langle\Lambda;\gamma\rangle = K\langle T_{\lambda}:\lambda\in
F\rangle /I K\langle T_{\lambda}:\lambda\in F\rangle \ .
$$
This means that the affinoid variety
$Sp(K\langle\Lambda;\gamma\rangle)$ is the preimage under $\iota$
of the affinoid unit polydisk in $\Adss^f$. In particular,
$Sp(K\langle\Lambda;\gamma\rangle)$ is an open subdomain in $T'$
which is reduced and coincides with the rational subset
$T'_\gamma$ (cf.\ [FvP] Prop.\ 4.6.1(4)). This establishes point
ii.\ of the assertion. The point iii. is Lemma 2.3. For point iv.,
as the invariants in an affinoid algebra with respect to a finite
group action, $K\langle \Lambda;\gamma\rangle^{W,\gamma}$ is again
affinoid (cf.\ [BGR] 6.3.3 Prop.\ 3).

\medskip

Suppose that the group $G$ is semisimple and adjoint. Then the
structure of the affinoid algebra
$K\langle\Lambda;\gamma\rangle^{W,\gamma}$ is rather simple. The
reason is that for such a group the set $\Lambda^{--}$ is the free
commutative monoid over the fundamental antidominant cocharacters
$\lambda_1,\dots,\lambda_d$. As usual we let $K\langle
X_1,\ldots,X_d\rangle$ denote the Tate algebra in $d$ variables
over $K$. Obviously we have a unique continuous algebra
homomorphism
$$
K\langle X_1,\ldots,X_d\rangle \longrightarrow
K\langle\Lambda;\gamma\rangle^{W,\gamma}
$$
sending the variable $X_i$ to $\sigma_{\lambda_i}$.

We also need a general lemma about orthogonal bases in normed
vector spaces. Let $(Y,\|\ \|)$ be a normed $K$-vector space and
suppose that $Y$ has an orthogonal basis of the form
$\{x_\ell\}_{\ell\in I}$. Recall that the latter means that
$$
\|\sum_\ell c_\ell x_\ell\| = \sup_\ell |c_\ell|\cdot\|x_\ell\|
$$
for any vector $\sum_\ell c_\ell x_\ell \in Y$. We suppose
moreover that there is given a partial order $\leq$ on the index
set $I$ such that:

-- Any nonempty subset of $I$ has a minimal element;

-- for any $k \in I$ the set $\{\ell \in I : \ell \leq k\}$ is
finite.

(Note that the partial order $\leq$ on $\Lambda^{--}$ has these
properties.)

\medskip

{\bf Lemma 2.5:} {\it Suppose that $\|x_\ell\| \leq \|x_k\|$
whenever $\ell \leq k$; furthermore, let elements $c_{\ell k} \in
K$ be given, for any $\ell \leq k$ in $I$, such that $|c_{\ell k}|
\leq 1$; then the vectors
$$
y_k := x_k + \sum_{\ell < k} c_{\ell k}x_\ell
$$
form another orthogonal basis of $Y$, and $\|y_k\| = \|x_k\|$.}

Proof: We have
$$
\|y_k\| = \max(\|x_k\|,\max_{\ell < k} |c_{\ell
k}|\cdot\|x_\ell\|) = \|x_k\|
$$
as an immediate consequence of our assumptions. We also have
$$
x_k = y_k + \sum_{\ell < k} b_{\ell k} y_\ell
$$
where $(b_{\ell k})$ is the matrix inverse to $(c_{\ell k})$ (over
the ring of integers in $K$; cf.\ [B-GAL] Chap. VI \S3.4 Lemma 4).
Let now $x = \sum_k c_k x_k$ be an arbitrary vector in $Y$. We
obtain
$$
x = \sum_k c_k x_k = \sum_k c_k (\sum_{\ell \leq k} b_{\ell k}
y_\ell) = \sum_\ell (\sum_{\ell \leq k} c_k b_{\ell k})y_\ell \ .
$$
Clearly $\|x\| \leq \sup_\ell |\sum_{\ell \leq k} c_k b_{\ell k}|
\cdot\|y_\ell\|$. On the other hand we compute
$$
\matrix{
 \sup_\ell |\sum_{\ell \leq k} c_k b_{\ell k}|
\cdot\|y_\ell\| & \leq\ \sup_\ell \sup_{\ell \leq k} |c_k|\cdot
\|y_\ell\| = \sup_\ell \sup_{\ell \leq k} |c_k|\cdot \|x_\ell\|
\cr & \leq\ \sup_k |c_k|\cdot \|x_k\| = \|x\| \ .\hfill }
$$

\medskip

{\bf Proposition 2.6:} {\it If the group $G$ is semisimple and
adjoint then the above map is an isometric isomorphism $K\langle
X_1,\ldots,X_d\rangle \mathop{\longrightarrow}\limits^{\cong}
K\langle\Lambda;\gamma\rangle^{W,\gamma}$.}

Proof: We write a given $\lambda \in \Lambda^{--}$ as $\lambda =
\lambda_1^{m_1}\ldots\lambda_d^{m_d}$ and put
$$
\widetilde{\sigma}_\lambda :=
\sigma_{\lambda_1}^{m_1}\cdot\ldots\cdot \sigma_{\lambda_d}^{m_d}
\ .
$$
It suffices to show that these
$\{\widetilde{\sigma}_\lambda\}_{\lambda \in \Lambda^{--}}$ form
another orthonormal basis of
$K\langle\Lambda;\gamma\rangle^{W,\gamma}$. One checks that the
arguments in [B-GAL] Chap. VI \S\S3.2 and 3.4 work, over the ring
of integers in $K$, equally well for our twisted $W$-action and
show that we have
$$
\widetilde{\sigma}_\lambda = \sigma_\lambda + \sum_{\mu < \lambda}
c_{\mu\lambda} \sigma_\mu
$$
with $|c_{\mu\lambda}| \leq 1$. So we may apply Lemma 2.5.

\medskip

We finish this section with a discussion of those examples of a
cocycle $\gamma$ which will be relevant later on.

\medskip

{\bf Example 1:} We fix a prime element $\pi_L$ of $L$. Let $\xi
\in X^\ast(T)$ be a dominant integral weight and put
$$
\gamma(w,\lambda(t)) := \pi_L^{\val_L(\xi({^w}t)) -
\val_L(\xi(t))} \ .
$$
This map $\gamma$ obviously has the properties $(a)$,$(b)$, and
$(d)$. For $t \in T^{--}$ we have $\lambda({^w}t) \leq \lambda(t)$
by [B-GAL] Chap.\ VI \S1.6 Prop.\ 18; since $\xi$ is dominant we
obtain
$$
\val_L \circ\xi ({t \over {^w}t}) \leq 0  \ .
$$
This means that $|\gamma(w,\lambda)| \leq 1$ for $\lambda \in
\Lambda^{--}$ which is condition $(c)$. We leave it as an exercise
to the reader to check that the resulting Banach algebra $K\langle
\Lambda;\gamma \rangle$ together with the twisted $W$-action, up
to isomorphism, is independent of the choice of the prime element
$\pi_L$.

{\bf Example 2:} A particular case of a dominant integral weight
is the determinant of the adjoint action of $T$ on the Lie algebra
$\Lie(N)$ of the unipotent radical $N$
$$
\Delta(t) := {\rm det}({\rm ad}(t);\Lie(N)) \ .
$$
Its absolute value satisfies
$$
\delta(t) = |\Delta(t)|_L^{-1}
$$
where $\delta : P \longrightarrow \Qdss^\times \subseteq K^\times$
is the modulus character of the Borel subgroup $P$. We let $K_q/K$
denote the splitting field of the polynomial $X^2 - q$ and we fix
a root $q^{1/2} \in K_q^\times$. Then the square root
$\delta^{1/2} : \Lambda \longrightarrow K_q^\times$ of the
character $\delta$ is well defined. For a completely analogous
reason as in the first example the cocycle
$$
\gamma(w,\lambda) := {\delta^{1/2}({^w}\lambda) \over
\delta^{1/2}(\lambda)}
$$
has the properties $(a) - (d)$. Moreover, using the root space
decomposition of $\Lie(N)$ one easily shows that
$$
\gamma(w,\lambda(t)) = \prod_{\alpha \in \Phi^+ \setminus
{^{w^{-1}}}\Phi^+} |\alpha(t)|_L \ .
$$
Hence the values of this cocycle $\gamma$ are integral powers of
$q$ and therefore lie in $K$.

{\bf Example 3:} Obviously the properties $(a) - (d)$ are preserved
by the product of two cocycles. For any dominant integral weight
$\xi \in X^\ast(T)$ therefore the cocycle
$$
\gamma_\xi(w,\lambda(t)) := {\delta^{1/2}({^w}\lambda) \over
\delta^{1/2}(\lambda)} \cdot \pi_L^{\val_L(\xi({^w}t)) -
\val_L(\xi(t))}
$$
is $K$-valued and satisfies $(a) - (d)$. We write
$$
V_\Rdss^\xi := V_\Rdss^{\gamma_\xi} \qquad\hbox{and}\qquad T'_\xi
:= T'_{\gamma_\xi} \ .
$$
Let $\eta \in V_\Rdss$ denote half the sum of the positive roots
in $\Phi^+$ and put
$$
\eta_L := [L:\Qdss_p]\cdot\eta \ .
$$
Let
$$
\xi_L := \val_L\circ\xi \in V_\Rdss \ .
$$
For the points $z_w \in V_\Rdss$ corresponding to the cocycle
$\gamma_\xi$ we then have
$$
z_w = (\eta_L + \xi_L) - {^w}(\eta_L + \xi_L) \ .
$$
In particular, $V_\Rdss^\xi$ is the convex hull of the points
${^w}(\eta_L + \xi_L) - (\eta_L + \xi_L)$ for $w \in W$. Note
that, since $\gamma_\xi$ has values in $L^\times$, the affinoid
variety $T'_\xi$ is naturally defined over $L$. Given any point $z
\in V_\Rdss$, we will write $z^{dom}$ for the unique dominant
point in the $W$-orbit of $z$.

\medskip

{\bf Lemma 2.7:} $V_\Rdss^\xi = \{ z \in V_\Rdss : (z + \eta_L +
\xi_L)^{dom} \leq \eta_L + \xi_L\}$.

Proof: Using the formula before Lemma 2.3 we have
$$
\eqalign{
 V_\Rdss^\xi & = \{ z \in V_\Rdss : {^w}z \leq (\eta_L +
 \xi_L) - {^w}(\eta_L + \xi_L)\ \hbox{for any}\ w \in W\}\cr
 & = \{ z \in V_\Rdss : {^w}(z + \eta_L + \xi_L) \leq \eta_L +
 \xi_L\ \hbox{for any}\ w \in W\}\ . }
$$
It remains to recall ([B-GAL] Chap.\ VI \S1.6 Prop.\ 18) that for
any $z \in V_\Rdss$ and any $w \in W$ one has ${^w}z \leq
z^{dom}$.

\medskip

The $\gamma_\xi$ in Example 3 are the cocycles which will appear
in our further investigation of specific Banach-Hecke algebras. In
the following we explicitly compute the affinoid domain $T'_\xi$
in case of the group $G := GL_{d+1}(L)$. (In case $\xi = 1$
compare also [Vig] Chap.\ 3.) We let $P \subseteq G$ be the lower
triangular Borel subgroup and $T \subseteq P$ be the torus of
diagonal matrices. We take $U_0 := GL_{d+1}(o_L)$ where $o_L$ is
the ring of integers of $L$. If $\pi_L \in o_L$ denotes a prime
element then
$$
\Lambda^{--} = \{\pmatrix{ \pi_L^{m_1}&&0 \cr & \ddots \cr
0&&\pi_L^{m_{d+1}} }T_0 : m_1 \geq \ldots \geq m_{d+1}\}\ .
$$
For $1 \leq i \leq d+1$ define the diagonal matrix
$$
t_i := \pmatrix{ \pi_L&&&&&0 \cr & \ddots \cr &&\pi_L \cr &&& 1
\cr &&&&\ddots \cr 0&&&&& 1 \cr}\quad\hbox{with $i$ diagonal
entries equal to $\pi_L$}\ .
$$
As a monoid $\Lambda^{--}$ is generated by the elements
$\lambda_1,\ldots,\lambda_{d+1},\lambda_{d+1}^{-1}$ where
$\lambda_i := \lambda(t_i)$. For any nonempty subset $I =
\{i_1,\ldots,i_s\} \subseteq \{1,\ldots,d+1\}$ let $\lambda_I \in
\Lambda$ be the cocharacter corresponding to the diagonal matrix
having $\pi_L$ at the places $i_1,\ldots,i_s$ and $1$ elsewhere.
Moreover let, as usual, $|I| := s$ be the cardinality of $I$ and
put $ht(I) := i_1 + \ldots + i_s$. These $\lambda_I$ together with
$\lambda_{\{1,\ldots,d+1\}}^{-1}$ form the $W$-orbit of the above
monoid generators. From the proof of Prop.\ 2.4 we therefore know
that $T'_\xi$ as a rational subdomain of $T'$ is described by the
conditions
$$
|\zeta(\lambda_I)| \leq |\gamma_\xi^{dom}(\lambda_I)|
$$
for any $I$ and
$$
|\zeta(\lambda_{\{1,\ldots,d+1\}})| =
|\gamma_\xi^{dom}(\lambda_{\{1,\ldots,d+1\}})| \ .
$$
One checks that
$$
|\gamma_1^{dom}(\lambda_I)| = |q|^{|I|(|I|+1)/2 - ht(I)} \ .
$$
If the dominant integral weight $\xi \in X^\ast(T)$ is given by
$$
\pmatrix{ g_1 &&0 \cr & \ddots \cr 0&& g_{d+1}} \longmapsto
\prod_{i=1}^{d+1} g_i^{a_i}
$$
with $(a_1,\ldots,a_{d+1}) \in \Zdss^{d+1}$ then
$$
|\gamma_\xi^{dom}(\lambda_I)| = |q|^{|I|(|I|+1)/2 - ht(I)}
 |\pi_L|^{\sum_{j=1}^{|I|} a_j - \sum_{i \in I} a_i} \ .
$$
We now use the coordinates
$$
\matrix{
 T'(K) & \longrightarrow & (K^\times)^{d+1} \hfill\cr
 \hfill \zeta & \longmapsto & (\zeta_1,\ldots,\zeta_{d+1}) \
 \hbox{with}\ \zeta_i := q^{i-1}\pi_L^{a_i}\zeta(\lambda_{\{i\}})
 }
$$
on the dual torus. In these coordinates $T'_\xi$ is the rational
subdomain of all $(\zeta_1,\ldots,\zeta_{d+1}) \in
(K^\times)^{d+1}$ such that
$$
\prod_{i \in I} |\zeta_i| \leq
|q|^{|I|(|I|-1)/2}|\pi_L|^{\sum_{i=1}^{|I|} a_i}
$$
for any proper nonempty subset $I \subseteq \{1,\ldots,d+1\}$ and
$$
\prod_{i=1}^{d+1} |\zeta_i| =
|q|^{d(d+1)/2}|\pi_L|^{\sum_{i=1}^{d+1} a_i} \ .
$$
The advantage of these variables is the following. As usual we
identify the Weyl group $W$ with the symmetric group on the set
$\{1,\ldots,d+1\}$. One checks that
$$
\gamma_\xi(w,\lambda_{\{i\}}) = q^{w(i)-i}\pi_L^{a_{w(i)} - a_i}
$$
for any $w \in W$ and $1 \leq i \leq d+1$. This implies that the
twisted $W$-action on the affinoid algebra $K\langle
\Lambda;\gamma_\xi \rangle$ is induced by the permutation action on
the coordinates $\zeta_1,\ldots,\zeta_{d+1}$ of the affinoid domain
$T'_\xi$. In fact, the above identity means that the cocycle
$\gamma_\xi$ can be written as the coboundary of an element in
$T'(K)$. This is more generally possible for any group $G$ whose
derived group is simply connected (cf.\ [Gro] \S8). We do not pursue
this point of view systematically, though, since it is not
compatible with general Langlands functoriality. But the problem of
``splitting'' the cocycle and the difficulty of reconciling the
normalization of the Satake isomorphism will reappear as a technical
complication in our attempt, in section 6, to treat Langlands
functoriality.

\medskip
\goodbreak

{\bf 3. The $p$-adic Satake isomorphism}

\medskip

Keeping the notations and assumptions introduced in the previous
section we now consider a locally $L$-analytic representation
$(\rho,E)$ of $G$ of the form
$$
E = K_\chi \otimes_L E_L
$$
where

-- $K_\chi$ is a one dimensional representation of $G$ given by a
locally $L$-analytic character $\chi : G \longrightarrow
K^\times$, and

-- $E_L$ is an $L$-rational irreducible representation $\rho_L$ of
$G$ of highest weight $\xi$.

Let
$$
E_L = \oplus_{\beta \in X^\ast(T)} E_{L,\beta}
$$
be the decomposition into weight spaces for $T$. According to [BT]
II.4.6.22 and Prop.\ II.4.6.28(ii) the reductive group $G$ has a
smooth connected affine model $\Gscr$ over the ring of integers
$o_L$ in $L$ such that $\Gscr(o_L) = U_0$. We fix once and for all
a $U_0$-invariant $o_L$-lattice $M$ in $E_L$ ([Jan] I.10.4) and
let $\|\ \|$ be the corresponding $U_0$-invariant norm on $E$. The
following fact is well-known.

\medskip

{\bf Lemma 3.1:} {\it We have $M =
\mathop{\oplus}\limits_{\beta\in X^\ast(T)} M_{\beta}$ with
$M_{\beta}:=M \cap E_{L,\beta}$.}

Proof: For the convenience of the reader we sketch the argument.
Fix a weight $\beta\in X^\ast(T)$. It suffices to construct an
element $\Pi_{\beta}$ in the algebra of distributions
$\Dist(\Gscr)$ which acts as a projector
$$
\Pi_{\beta}:E_L \longrightarrow E_{L,\beta} \ .
$$
Let $B$ be the finite set of weights $\neq\beta$ which occur in
$E_L$. Also we need the Lie algebra elements
$$
H_i:=(d\mu_i)(1)\in\Lie(\Gscr)
$$
where $\mu_1,\ldots,\mu_r$ is a basis of the cocharacter group of
$T$. We have
$$
\underline{\gamma}:=(d\gamma(H_1),\ldots,d\gamma(H_r))\in\Z^r\ \
\hbox{for any}\ \gamma\in X^\ast(T) \ .
$$
According to [Hum] Lemma 27.1 we therefore find a polynomial
$\Pi\in\Q [ y_1{,}{\ldots}{,}y_r]$ such that
$\Pi(\Z^r)\subseteq\Z$, $\Pi(\underline{\beta})=1$, and
$\Pi(\underline{\gamma})=0$ for any $\gamma\in B$. Moreover [Hum]
Lemma 26.1 says that the polynomial $\Pi$ is a $\Z$-linear
combination of polynomials of the form
$$
\pmatrix{ y_1\cr b_1\cr} \cdot\ldots\cdot \pmatrix{ y_r\cr b_r\cr}
\ \ \hbox{with integers}\ b_1,\ldots, b_r\geq 0 \ .
$$
Then [Jan] II.1.12 implies that
$$
\Pi_{\beta}:=\Pi(H_1,\ldots,H_r)
$$
lies in $\Dist(\Gscr)$. By construction $\Pi_{\beta}$ induces a
projector from $E_L$ onto $E_{L,\beta}$.

\medskip

It follows that, for any $t \in T$, the operator norm of
$\rho_L(t)$ on $E_L$ is equal to
$$
\|\rho_L(t)\| = \max\{|\beta(t)| : \beta \in X^\ast(T)\ \hbox{such
that}\ E_{L,\beta} \neq 0\}.
$$

\medskip

{\bf Lemma 3.2:} {\it For any $t \in T$ we have $\|\rho(t)\| =
|\chi(t)| \cdot |\xi({^w}t)|$ with $w \in W$ such that ${^w}t \in
T^{--}$.}

Proof: Consider first the case $t \in T^{--}$ with $w = 1$. For
any weight $\beta$ occurring in $E_L$ one has $\xi = \alpha\beta$
where $\alpha$ is an appropriate product of simple roots. But by
definition of $T^{--}$ we have $|\alpha(t)|_L \geq 1$ for any
simple root $\alpha$. For general $t \in T$ and $w \in W$ as in
the assertion we then obtain
$$
\eqalign{
 |\xi({^w}t)| & = \max\{|\beta({^w}t)| : E_{L,\beta} \neq 0\} \cr
              & = \max\{|\beta(t)| : E_{L,\beta} \neq 0\} \cr
              & = \|\rho_L(t)\| \ . }
$$
Here the second identity is a consequence of the fact that the set
of weights of $E_L$ is $W$-invariant.

\medskip

Collecting this information we first of all see that Lemma 1.4
applies and gives, for any open subgroup $U \subseteq U_0$, the
isomorphism
$$
\Hscr(G,1_U) \cong \Hscr(G,\rho|U) \ .
$$
But the norm $\|\ \|$ on $\Hscr(G,\rho|U)$ corresponds under this
isomorphism to the norm $\|\ \|_{\chi,\xi}$ on $\Hscr(G,1_U)$
defined by
$$
\|\psi\|_{\chi,\xi} := \sup_{g \in G}
|\psi(g)\chi(g)|\cdot\|\rho_L(g)\| \ .
$$
If $|\chi|$= 1 (e.g., if the group $G$ is semisimple) then the
character $\chi$ does not affect the norm $\|\ \|_\xi := \|\
\|_{\chi,\xi}$. In general $\chi$ can be written as a product
$\chi = \chi_1 \chi_{un}$ of two characters where $|\chi_1| = 1$
and $\chi_{un}|U_0 = 1$. Then
$$
\matrix{
 (\Hscr(G,1_U),\|\ \|_\xi) &
 \mathop{\longrightarrow}\limits^{\cong} & (\Hscr(G,1_U),\|\ \|_{\chi,\xi}) \cr
 \hfill \psi & \longmapsto & \psi\cdot\chi_{un}^{-1} \hfill }
$$
is an isometric isomorphism. We therefore have the following
fact.

\medskip

{\bf Lemma 3.3:} {\it The map
$$
\matrix{
 \|\ \|_\xi\hbox{-completion of}\ \Hscr(G,1_U) &
 \mathop{\longrightarrow}\limits^{\cong} &
 \Bscr(G,\rho|U) \cr
 \hfill \psi & \longmapsto & \psi\cdot \chi_{un}^{-1}\rho \hfill }
$$
is an isometric isomorphism of Banach algebras.}

\medskip

In this section we want to compute these Banach-Hecke algebras in
the case $U = U_0$. By the Cartan decomposition $G$ is the
disjoint union of the double cosets $U_0t U_0$ with $t$ running
over $T^{--}/T_0$. Let therefore $\psi_{\lambda(t)} \in
\Hscr(G,1_{U_0})$ denote the characteristic function of the double
coset $U_0tU_0$. Then $\{\psi_\lambda\}_{\lambda \in
\Lambda^{--}}$ is a $K$-basis of $\Hscr(G,1_{U_0})$. According to
Lemma 3.2 the norm $\|\ \|_\xi$ on $\Hscr(G,1_{U_0})$ is given by
$$
\|\psi\|_\xi := \sup_{t \in T^{--}} |\psi(t)\xi(t)| \ .
$$
The $\{\psi_\lambda\}_{\lambda \in \Lambda^{--}}$ form a $\|\
\|_\xi$-orthogonal basis of $\Hscr(G,1_{U_0})$ and hence of its
$\|\ \|_\xi$-completion.

The Satake isomorphism computes the Hecke algebra
$\Hscr(G,1_{U_0})$. For our purposes it is important to consider
the renormalized version of the Satake map given by
$$
\matrix{
 S_\xi : \Hscr(G,1_{U_0}) & \longrightarrow & K[\Lambda] \hfill\cr
 \hfill \psi & \longmapsto & \mathop{\sum}\limits_{t \in T/T_0}
 \pi_L^{\val_L(\xi(t))}(\mathop{\sum}\limits_{n \in
 N/N_0} \psi(tn))\lambda(t) \ .}
$$
On the other hand we again let $K_q/K$ be the splitting field of
the polynomial $X^2 - q$ and we temporarily fix a root $q^{1/2}
\in K_q$. Satake's theorem says (cf. [Car]\S4.2) that the map
$$
\matrix{
 S^{norm} : \Hscr(G,1_{U_0}) \otimes_K K_q & \longrightarrow & K_q[\Lambda] \hfill\cr
 \hfill \psi & \longmapsto & \mathop{\sum}\limits_{t \in T/T_0}
 \delta^{-1/2}(t)(\mathop{\sum}\limits_{n \in
 N/N_0} \psi(tn))\lambda(t) }
$$
induces an isomorphism of $K_q$-algebras
$$
\Hscr(G,1_{U_0}) \otimes_K K_q
\mathop{\longrightarrow}\limits^{\cong} K_q[\Lambda]^W \ .
$$
Here the $W$-invariants on the group ring $K_q[\Lambda]$ are
formed with respect to the $W$-action induced by the conjugation
action of $N(T)$ on $T$. Since
$\pi_L^{\val_L\circ\xi}\delta^{1/2}$ defines a character of
$\Lambda$ it is clear that $S_\xi$ is a homomorphism of algebras
as well and a simple Galois descent argument shows that $S_\xi$
induces an isomorphism of $K$-algebras
$$
\Hscr(G,1_{U_0}) \mathop{\longrightarrow}\limits^{\cong}
K[\Lambda]^{W,\gamma_\xi}
$$
where $\gamma_\xi$ is the cocycle from Example 3 in section 2. The
left hand side has the $\|\ \|_\xi$-orthogonal basis
$\{\psi_\lambda\}_{\lambda \in \Lambda^{--}}$ with
$$
\|\psi_{\lambda(t)}\|_\xi = |\xi(t)| \ .
$$
The right hand side has the $\|\ \|_{\gamma_\xi}$-orthonormal
basis $\{\sigma_\lambda\}_{\lambda \in \Lambda^{--}}$ where
$$
\sigma_\lambda = \sum_{w \in W/W(\lambda)}
\gamma_\xi(w,\lambda){^w}\lambda
$$
(cf.\ section 2). Since the maps
$$
\matrix{
 N/N_0 & \mathop{\longrightarrow}\limits^{\simeq} & NtU_0/U_0 \cr
 \hfill nN_0 & \longmapsto & tnU_0 \hfill }
$$ are bijections we have
$$
\sum_{n \in N/N_0} \psi_{\lambda(s)}(tn) = |(NtU_0 \cap
U_0sU_0)/U_0| =: c(\lambda(t),\lambda(s))\qquad\hbox{for any}\ s, t
\in T\ .
$$
It follows that
$$
\eqalign{
 S_\xi(\psi_\mu) & =
    \mathop{\sum}\limits_{t \in T/T_0}
    \pi_L^{\val_L(\xi(t))}c(\lambda(t),\mu)\lambda(t) \cr
    & =
    \mathop{\sum}\limits_{\lambda \in \Lambda^{--}}
    \pi_L^{\val_L\circ\xi(\lambda)}c(\lambda,\mu)\sigma_\lambda
    \qquad\hbox{for any}\ \mu \in
    \Lambda^{--}\ . }
$$
and
$$
\pi_L^{\val_L\circ\xi({^w}\lambda)}c({^w}\lambda,\mu) =
\gamma_\xi(w,\lambda)\pi_L^{\val_L\circ\xi(\lambda)}c(\lambda,\mu)
$$
for any $\lambda \in \Lambda^{--}$, $\mu \in \Lambda$, and $w \in
W$.

The reason for the validity of Satake's theorem lies in the
following properties of the coefficients $c(\lambda,\mu)$.

\medskip

{\bf Lemma 3.4:} {\it For $\lambda, \mu \in \Lambda^{--}$ we have:

i. $c(\mu,\mu) = 1$;

ii. $c(\lambda,\mu) = 0$ unless $\lambda \leq \mu$.}

Proof: [BT] Prop.\ I.4.4.4.

\medskip

{\bf Proposition 3.5:} {\it The map $S_\xi$ extends by continuity
to an isometric isomorphism of $K$-Banach algebras}
$$
\|\ \|_\xi\hbox{\it -completion of}\ \Hscr(G,1_{U_0})
\mathop{\longrightarrow}\limits^{\cong}
K\langle\Lambda;\gamma_\xi\rangle^{W,\gamma_\xi}\ .
$$

Proof: Define
$$
\widetilde{\psi}_{\lambda} :=
\pi_L^{-\val_L\circ\xi(\lambda)}\psi_{\lambda}
$$
for $\lambda \in \Lambda^{--}$. The left, resp. right, hand side
has the $\|\ \|_\xi$-orthonormal, resp. $\|\
\|_{\gamma_\xi}$-orthonormal, basis
$\{\widetilde{\psi}_\lambda\}_{\lambda\in\Lambda^{--}}$, resp.
$\{\sigma_\lambda\}_{\lambda\in\Lambda^{--}}$. We want to apply
Lemma 2.5 to the normed vector space
$(K[\Lambda]^{W,\gamma_\xi},\|\ \|_{\gamma_\xi})$, its orthonormal
basis $\{\sigma_\lambda\}$, and the elements
$$
S_\xi(\widetilde{\psi}_ \mu) = \sigma_\mu + \sum_{\lambda < \mu}
\pi_L^{\val_L\circ\xi(\lambda) -
\val_L\circ\xi(\mu)}c(\lambda,\mu)\sigma_\lambda
$$
(cf. Lemma 3.4). The coefficients $c(\lambda,\mu)$ are integers
and therefore satisfy $|c(\lambda,\mu)| \leq 1$. Moreover,
$\lambda < \mu$ implies, since $\xi$ is dominant, that
$\val_L\circ\xi(\mu) \leq \val_L\circ\xi(\lambda)$. Hence the
assumptions of Lemma 2.5 indeed are satisfied and we obtain that
$\{S_\xi(\widetilde{\psi}_\lambda)\}$ is another orthonormal basis
for $(K[\Lambda]^{W,\gamma_\xi},\|\ \|_{\gamma_\xi})$.

\medskip

{\bf Corollary 3.6:} {\it The Banach algebras $\Bscr(G,\rho|U_0)$
and $K\langle\Lambda;\gamma_\xi\rangle^{W,\gamma_\xi}$ are
isometrically isomorphic.}

\medskip

If $\xi = 1$ then, in view of Lemma 2.7, the reader should note
the striking analogy between the above result and the computation
in [Mac] Thm.\ (4.7.1) of the spectrum of the algebra of
integrable complex valued functions on $U_0\backslash G/U_0$. The
methods of proof are totally different, though. In fact, in our
case the spherical function on $U_0\backslash G/U_0$ corresponding
to a point in $T'_1$ in general is not bounded.

Suppose that the group $G$ is semisimple and adjoint. We fix
elements $t_1,\ldots,t_d \in T^{--}$ such that $\lambda_i :=
\lambda(t_i)$ are the fundamental antidominant cocharacters. In
Prop.\ 2.6 we have seen that then
$K\langle\Lambda;\gamma_\xi\rangle^{W,\gamma_\xi}$ is a Tate
algebra in the variables
$\sigma_{\lambda_1}\ldots,\sigma_{\lambda_d}$. Hence
$\Bscr(G,\rho|U_0)$ is a Tate algebra as well.  But it seems
complicated to compute explicitly the variables corresponding to
the $\sigma_{\lambda_i}$. Instead we may repeat our previous
reasoning in a modified way.

\medskip

{\bf Proposition 3.7:} {\it Suppose that $G$ is semisimple and
adjoint; then $\Bscr(G,\rho|U_0)$ is a Tate algebra over $K$ in
the variables ${\psi_{\lambda_1}\cdot\rho \over \xi(t_1)},\ldots,
{\psi_{\lambda_d}\cdot\rho \over \xi(t_d)}$.}

Proof: By Lemma 3.3 and Prop.\ 3.5 it suffices to show that
$K\langle\Lambda;\gamma_\xi\rangle^{W,\gamma_\xi}$ is a Tate
algebra in the variables $\xi(t_i)^{-1}S_\xi(\psi_{\lambda_i})$.
We write a given $\lambda \in \Lambda^{--}$ as $\lambda =
\lambda_1^{m_1}\ldots\lambda_d^{m_d}$ and put
$$
\widetilde{\sigma}_\lambda :=
S_\xi(\widetilde{\psi}_{\lambda_1})^{m_1}\cdot\ldots\cdot
S_\xi(\widetilde{\psi}_{\lambda_d})^{m_d} =
S_\xi(\widetilde{\psi}_{\lambda_1}^{m_1}\ast\ldots\ast
\widetilde{\psi}_{\lambda_d}^{m_d})
$$
using notation from the proof of Prop.\ 3.5. Similarly as in the
proof of Prop.\ 2.7 one checks that the arguments in [B-GAL] Chap.
VI \S\S3.2 and 3.4 work, over the ring of integers in $K$, equally
well for our twisted $W$-action (note that, in the language of
loc.\ cit.\ and due to Lemma 3.4, the unique maximal term in
$S_\xi(\widetilde{\psi}_{\lambda_i})$ is $\lambda_i$) and show
that we have
$$
\widetilde{\sigma}_\lambda = \sigma_\lambda + \sum_{\mu < \lambda}
c_{\mu\lambda} \sigma_\mu
$$
with $|c_{\mu\lambda}| \leq 1$. So we may apply again Lemma 2.5
and obtain that $\{\widetilde{\sigma}_\lambda\}$ is another
orthonormal basis for
$K\langle\Lambda;\gamma_\xi\rangle^{W,\gamma_\xi}$. It remains to
note that $\xi(t_i)$ and $\pi_L^{\val_L(\xi(t))}$ only differ by a
unit.

\medskip

{\bf Example:} Consider the group $G := GL_{d+1}(L)$. Cor.\ 3.6
applies to $G$ but Prop.\ 3.7 does not. Nevertheless, with the
same notations as at the end of section 2 a simple modification of
the argument gives
$$
\Bscr(G,\rho|U_0) =
K\big\langle{\psi_{\lambda_1}\cdot\chi_{un}^{-1}\rho \over
\xi(t_1)},\ldots, {\psi_{\lambda_d}\cdot\chi_{un}^{-1}\rho \over
\xi(t_d)}, \big({\psi_{\lambda_{d+1}}\cdot\chi_{un}^{-1}\rho \over
\xi(t_{d+1})}\big)^{\pm 1} \big\rangle \ .
$$
Moreover in this case the $\lambda_i$ are minimal with respect to
the partial order $\leq$ so that we do have
$$
\xi(t_i)^{-1}S_\xi(\psi_{\lambda_i}) = \sigma_{\lambda_i} \ .
$$
Hence the above representation of $\Bscr(G,\rho|U_0)$ as an
affinoid algebra corresponds to the representation
$$
K\langle\Lambda;\gamma_\xi\rangle^{W,\gamma_\xi} = K\langle
\sigma_{\lambda_1}\ldots,\sigma_{\lambda_d},\sigma_{\lambda_{d+1}}^{\pm
1} \rangle \ .
$$
On affinoid domains this corresponds to a map
$$
T'_\xi \longrightarrow \{(\omega_1,\ldots,\omega_{d+1}) \in
K^{d+1} : |\omega_1|,\ldots,|\omega_d| \leq 1, |\omega_{d+1}| =
1\}
$$
which, using our choice of coordinates on $T'$ from section 2, is
given by
$$
(\zeta_1,\ldots,\zeta_{d+1}) \longmapsto (\ldots,q^{(i-1)i \over
2}\xi(t_i)^{-1}\Sigma_i(\zeta_1,\ldots,\zeta_{d+1}),\ldots)
$$
where
$$
\Sigma_1(\zeta_1,\ldots,\zeta_{d+1}) = \zeta_1 + \ldots +
\zeta_{d+1},\ \ldots\ , \Sigma_{d+1}(\zeta_1,\ldots,\zeta_{d+1}) =
\zeta_1 \cdot\ldots\cdot \zeta_{d+1}
$$
denote the elementary symmetric polynomials.

Let us further specialize to the case $G = GL_2(L)$. Then $E_L$ is
the $k$-th symmetric power, for some $k \geq 0$, of the standard
representation of $GL_2$. The highest weight of $E_L$ is
$\xi(\pmatrix{t_1 & 0\cr 0 & t_2}) = t_2^k$. We obtain
$$
\Bscr(G,\rho|U_0) = K\big\langle X_1, (\pi_L^{-k}X_2)^{\pm 1}
\big\rangle \ .
$$
with the variables $X_i :=
\psi_{\lambda_i}\cdot\chi_{un}^{-1}\rho$. The above map between
affinoid domains becomes
$$
(\zeta_1,\zeta_2) \longrightarrow (\zeta_1 + \zeta_2,
q^{-1}\pi_L^{-k} \zeta_1\zeta_2) \ .
$$

\medskip

{\bf 4. $p$-adic Iwahori-Hecke algebras}

\medskip

With the same assumptions and notations as in the previous section
we now let $U_1 \subseteq U_0$ be the Iwahori subgroup such that
$U_1 \cap P = U_0 \cap P$. In this section we will compute the
Banach-Hecke algebras $\Bscr(G,\rho|U_1)$. By Lemma 3.3 this
means, similarly as before, computing the $\|\ \|_\xi$-completion
of $\Hscr(G,1_{U_1})$.

The extended affine Weyl group $\widetilde{W}$ of $G$ is given by
$$
\widetilde{W} := N(T)/T_0 \ .
$$
Since the Weyl group $W$ lifts to $U_0 \cap N(T)/T_0 \subseteq
\widetilde{W}$ we see that $\widetilde{W}$ is the semidirect
product of $W$ and $\Lambda$. The Bruhat-Tits decomposition says
that $G$ is the disjoint union of the double cosets $U_1xU_1$ with
$x$ running over $\widetilde{W}$. Therefore, if we let $\tau_x \in
\Hscr(G,1_{U_1})$ denote the characteristic function of the double
coset $U_1xU_1$, then $\{\tau_x\}_{x \in \widetilde{W}}$ is a
$K$-basis of $\Hscr(G,1_{U_1})$. The $\tau_x$ are known to be
inver- tible in the algebra $\Hscr(G,1_{U_1})$. As a consequence
of Lemma 3.2 the $\|\ \|_\xi$-norm is given by
$$
\|\psi\|_\xi = \sup_{v,w \in W} \sup_{t \in T^{--}}
|\psi(v\lambda({^w}t))\xi(t)| \ .
$$
In particular, $\{\tau_x\}_{x \in \widetilde{W}}$ is an $\|\
\|_\xi$-orthogonal basis of $\Hscr(G,1_{U_1})$ such that
$$
\|\tau_x \|_\xi = |\xi({^w}t)| \qquad\hbox{if}\ v,w \in W\
\hbox{and}\ t \in T\ \hbox{such that}\ x = v\lambda(t)\
\hbox{and}\ {^w}t \in T^{--} \ .
$$
We let $\Cscr$ be the unique Weyl chamber corresponding to $P$ in
the apartment corresponding to $T$ with vertex $x_0$ (cf.\
[Car]\S3.5). The Iwahori subgroup $U_1$ fixes pointwise the unique
chamber $C \subseteq \Cscr$ with vertex $x_0$. The reflections at
the walls of $\Cscr$ generate the Weyl group $W$. Let
$s_0,\ldots,s_e \in \widetilde{W}$ be the reflections at all the
walls of $C$ and let $W_{aff}$ denote the subgroup of
$\widetilde{W}$ generated by $s_0,\ldots,s_e$. This affine Weyl
group $W_{aff}$ with the generating set $\{s_0,\ldots,s_e\}$ is a
Coxeter group. In particular we have the corresponding length
function $\ell : W_{aff} \longrightarrow \Ndss \cup \{0\}$ and the
corresponding Bruhat order $\leq$ on $W_{aff}$. If $\Omega
\subseteq \widetilde{W}$ is the subgroup which fixes the chamber
$C$ then $\widetilde{W}$ also is the semidirect product of
$\Omega$ and $W_{aff}$. We extend the length function $\ell$ to
$\widetilde{W}$ by $\ell(\omega w) := \ell(w)$ for $\omega \in
\Omega$ and $w \in W_{aff}$. The Bruhat order is extended to
$\widetilde{W}$ by the rule $\omega w \leq \omega' w'$, for $w,w'
\in W_{aff}$ and $\omega,\omega' \in \Omega$, if and only if
$\omega = \omega'$ and $w \leq w'$. One of the basic relations
established by Iwahori-Matsumoto is:

$(1)\quad$ For any $x,y \in \widetilde{W}$ such that $\ell(xy) =
\ell(x) + \ell(y)$ we have $\tau_{xy} = \tau_x \ast \tau_y$.

It easily implies that, for any $\lambda \in \Lambda$, the element
$$
\Theta(\lambda) := \tau_{\lambda_1}\ast\tau^{-1}_{\lambda_2} \in
\Hscr(G,1_{U_1})
$$
where $\lambda = \lambda_1 \lambda_2^{-1}$ with $\lambda_i \in
\Lambda^{--}$ is independent of the choice of $\lambda_1$ and
$\lambda_2$. Moreover Bernstein has shown that the map
$$
\matrix{
 \Theta : K[\Lambda] & \longrightarrow & \Hscr(G,1_{U_1}) \cr
 \hfill \lambda & \longmapsto & \Theta(\lambda) \hfill }
$$
is an embedding of $K$-algebras.

{\it Comment:} It is more traditional (cf.\ [HKP] \S1) to consider
the embedding of $K_q$-algebras (with $K_q/K$ and $q^{1/2} \in
K_q$ be as before)
$$
\matrix{
 \Theta^{norm} : K_q[\Lambda] & \longrightarrow & \Hscr(G,1_{U_1}) \otimes_K K_q \hfill\cr
 \hfill \lambda & \longmapsto &
 \delta^{-1/2}(\lambda)\tau_{\lambda_1}\ast\tau^{-1}_{\lambda_2}  \hfill }
$$
where $\lambda = \lambda_1 \lambda_2^{-1}$ with dominant
$\lambda_i$. The modified map $\Theta^+ :=
\delta^{1/2}\cdot\Theta^{norm}$ already is defined over $K$. On
$K[\Lambda]$ we have the involution $\iota_\lambda$ defined by
$\iota_\Lambda(\lambda) := \lambda^{-1}$, and on
$\Hscr(G,1_{U_1})$ there is the anti-involution $\iota$ defined by
$\iota(\psi)(g) := \psi(g^{-1})$. We then have
$$
\Theta = \iota \circ \Theta^+ \circ \iota_\Lambda \ .
$$

In the following we consider the renormalized embedding of
$K$-algebras
$$
\matrix{
 \Theta_\xi :\qquad K[\Lambda] & \longrightarrow & \Hscr(G,1_{U_1}) \hfill \cr
 \hfill \lambda & \longmapsto & \pi_L^{-\val_L\circ\xi(\lambda)}
 \Theta(\lambda) \ . \hfill }
$$
In order to compute the norm induced, via $\Theta_\xi$, by $\|\
\|_\xi$ on $K[\Lambda]$ we introduce the elements
$$
\theta_x :=
q^{(\ell(x)-\ell(w)-\ell(\lambda_1)+\ell(\lambda_2))/2} \tau_w
\ast \tau_{\lambda_1} \ast \tau_{\lambda_2}^{-1} \ .
$$
for any $x \in \widetilde{W}$ written as $x =
w\lambda_1\lambda_2^{-1}$ with $w \in W$ and $\lambda_i \in
\Lambda^{--}$. Since $\ell(w) + \ell(\lambda_1) =
\ell(w\lambda_1)$ (cf.\ [Vig] App.) we obtain from $(1)$ that
$$
\theta_x =
q^{(\ell(w\lambda_1\lambda_2^{-1})-\ell(w\lambda_1)+\ell(\lambda_2))/2}
\tau_{w\lambda_1} \ast \tau_{\lambda_2}^{-1} \ .
$$
On the other hand [Vig] Lemma 1.2 (compare also [Hai] Prop.\ 5.4)
says that, for any $x,y \in \widetilde{W}$, the number
$$
(\ell(xy^{-1}) - \ell(x) + \ell(y))/2
$$
is an integer between $0$ and $\ell(y)$ and that
$$
\tau_x \ast \tau_y^{-1} = q^{-(\ell(xy^{-1}) - \ell(x) +
\ell(y))/2} (\tau_{xy^{-1}} + Q_{x,y})
$$
where $Q_{x,y}$ is a linear combination with integer coefficients
of $\tau_z$ with $z < xy^{-1}$. It follows that for any $x \in
\widetilde{W}$ we have
$$
\theta_x = \tau_x + Q_x \leqno{(2)}
$$
where $Q_x$ is a linear combination with integer coefficients of
$\tau_z$ with $z < x$.

\medskip

{\bf Lemma 4.1:} {\it Consider two elements $x = w'\lambda$ and $y
= v'\mu$ in $\widetilde{W}$ where $w',v' \in W$ and $\lambda, \mu
\in \Lambda$; let $w,v \in W$ such that ${^w}\lambda , {^v}\mu \in
\Lambda^{--}$; if $x \leq y$ then we have:

i. ${^v}\mu - {^w}\lambda \in \sum_{\alpha \in \Phi^+} \Ndss_0
\cdot (-\check{\alpha})$;

ii. $\|\tau_x\|_\xi \leq \|\tau_y\|_\xi$.}

Proof: i. Let $w_0 \in W$ denote the longest element. We will make
use of the identity
$$
\{x' \in \widetilde{W} : x' \leq w_0({^{w_0v}}\mu)\} =
\bigcup_{\lambda'} W\lambda' W
$$
where $\lambda'$ ranges over all elements in $\Lambda^{--}$ such
that ${^v}\mu - \lambda' \in \sum_{\alpha \in \Phi^+} \Ndss_0 \cdot
(-\check{\alpha})$ (see [Ka2] (4.6) or [HKP] 7.8). Since $y \in W
({^v}\mu) W$ this identity implies first that $x \leq y \leq
({^{w_0v}}\mu) w_0$ and then that $x \in W\lambda' W$ for some
$\lambda' \in \Lambda^{--}$ such that ${^v}\mu - \lambda' \in
\sum_{\alpha \in \Phi^+} \Ndss_0 \cdot (-\check{\alpha})$. Obviously
we must have $\lambda' = {^v}\lambda$. ii. Let $\lambda =
\lambda(t_1)$ and $\mu = \lambda(t_2)$. We have $\|\tau_x\|_\xi
=|\xi({^w}t_1)|$ and $\|\tau_y\|_\xi = |\xi({^v}t_2)|$. Since
highest weights are dominant we obtain from i. that
$|\xi({^v}t_2({^w}t_1)^{-1})| \geq 1$.

\medskip

It follows from Lemma 4.1.ii and formula $(2)$ that Lemma 2.5 is
applicable showing that $\{\theta_x\}_{x \in \widetilde{W}}$ is
another $\|\ \|_\xi$-orthogonal basis of $\Hscr(G,1_{U_1})$ with
$$
\|\theta_x\|_\xi = \|\tau_x\|_\xi \ .
$$
For any $\lambda(t) = \lambda = \lambda_1\lambda_2^{-1} \in
\Lambda$ with $\lambda_i \in \Lambda^{--}$ we have
$$
\theta_\lambda =
q^{(\ell(\lambda)-\ell(\lambda_1)+\ell(\lambda_2))/2}\pi_L^{\val_L(\xi(t))}
\Theta_\xi(\lambda)
$$
and
$$
\|\theta_\lambda\|_\xi = \|\tau_\lambda \|_\xi = |\xi({^w}t)|
$$
where $w \in W$ such that ${^w}t \in T^{--}$. In particular
$\{\theta_\lambda\}_{\lambda \in \Lambda}$ is a $\|\
\|_\xi$-orthogonal basis of ${\rm im}(\Theta_\xi)$.

\medskip

{\bf Lemma 4.2:} {\it With the above notations we have}
$$
q^{-(\ell(\lambda)-\ell(\lambda_1)+\ell(\lambda_2))/2} =
{\delta^{1/2}({^w}\lambda) \over \delta^{1/2}(\lambda)} \ .
$$

Proof: Write $t = t_1 t_2^{-1}$ with $\lambda(t_i) = \lambda_i$.
According to the explicit formula for the length $\ell$ in [Vig]
App.\  we have
$$
q^{\ell(\lambda)} = \prod_{\alpha \in \Phi^+, |\alpha(t)|_L \geq
1} |\alpha(t)|_L \cdot \prod_{\alpha \in \Phi^+, |\alpha(t)|_L
\leq 1} |\alpha(t)|_L^{-1}
$$
and
$$
q^{\ell(\lambda_i)} = \prod_{\alpha \in \Phi^+} |\alpha(t_i)|_L \
.
$$
It follows that
$$
q^{-(\ell(\lambda)-\ell(\lambda_1)+\ell(\lambda_2))/2} =
\prod_{\alpha \in \Phi^+, |\alpha(t)|_L \leq 1} |\alpha(t)|_L
$$
Since ${^w}t \in T^{--}$ we have $|{^{w^{-1}}}\alpha(t)|_L \geq 1$
for any $\alpha \in \Phi^+$. Hence $\{\alpha \in \Phi^+ :
|\alpha(t)|_L < 1\} \subseteq \Phi^+ \setminus {^{w^{-1}}}\Phi^+$.
By the last formula in Example 2 of section 2 the above right hand
side therefore is equal to ${\delta^{1/2}({^w}\lambda) \over
\delta^{1/2}(\lambda)}$.

\medskip

It readily follows that
$$
\|\Theta_\xi(\lambda)\|_\xi = |\gamma_\xi^{dom}(\lambda)|
\qquad\hbox{for any}\ \lambda \in \Lambda \ .
$$
In other words
$$
 \Theta_\xi : (K[\Lambda],\|\ \|_{\gamma_\xi})  \longrightarrow
  (\Hscr(G,1_{U_1}),\|\ \|_\xi)
$$
is an isometric embedding. Combining all this with Lemma 3.3 we
obtain the following result.

\medskip

{\bf Proposition 4.3:} {\it i. The map
$$
\matrix{
 K\langle \Lambda;\gamma_\xi \rangle &
 \longrightarrow & \Bscr(G,\rho|U_1) \hfill\cr
 \hfill \lambda & \longmapsto & \Theta_\xi(\lambda)\cdot \chi_{un}^{-1}\rho }
$$
is an isometric embedding of Banach algebras;

ii. the map
$$
\matrix{
 \Hscr(U_0,1_{U_1}) \otimes_K K\langle \Lambda;\gamma_\xi \rangle &
 \mathop{\longrightarrow}\limits^{\cong} & \Bscr(G,\rho|U_1) \hfill\cr
 \hfill \tau_w \otimes \lambda & \longmapsto &
 (\tau_w \ast \Theta_\xi(\lambda))\cdot \chi_{un}^{-1}\rho }
$$
is a $K$-linear isomorphism.}

\medskip

{\bf Remarks:} 1) A related computation in the case $\xi = 1$ is
contained in [Vig] Thm.\ 4(suite).

2) It is worth observing that the ``twisted" $W$-action on
$K\langle\Lambda;\gamma_\xi\rangle$ corresponds under the
isomorphism $\Theta_\xi$ to the $W$-action on ${\rm
im}(\Theta_\xi)$ given by
$$
(w,\theta_{\lambda}) \longmapsto \theta_{^{w}\lambda} \ .
$$

\medskip

The results of this section and of the previous section are
compatible in the following sense.

\medskip

{\bf Proposition 4.4:} {\it The diagram
$$
\xymatrix{
  K\langle \Lambda;\gamma_\xi \rangle \ar[rr]^{\Theta_\xi(.)\cdot \chi_{un}^{-1}\rho\quad}
  & & \Bscr(G,\rho|U_1)
  \ar[d]^{(\psi_{\lambda(1)}\cdot \chi_{un}^{-1}\rho)\ast .} \\
  K\langle \Lambda;\gamma_\xi \rangle^{W,\gamma_\xi} \ar[u]^{\subseteq}
  \ar[rr]^{S_\xi^{-1}(.)\cdot \chi_{un}^{-1}\rho} & &
  \Bscr(G,\rho|U_0)   }
$$
is commutative. Moreover, the image of $K\langle
\Lambda;\gamma_\xi \rangle^{W,\gamma_\xi}$ under the map
$\Theta_\xi(.)\cdot \chi_{un}^{-1}\rho$ lies in the center of
$\Bscr(G,\rho|U_1)$.}

Proof: We recall that the upper, resp.\ lower, horizontal arrow is
an isometric unital monomorphism by Prop.\ 4.3.i, resp. by Lemma
3.3 and Prop. 3.5. The right perpendicular arrow is a continuous
linear map respecting the unit elements. It suffices to treat the
case of the trivial representation $\rho = 1$. By continuity we
therefore are reduced to establishing the commutativity of the
diagram
$$
\xymatrix{
  K[\Lambda] \ar[r]^{\Theta\quad} & \Hscr(G,1_{U_1})
  \ar[d]^{\psi_{\lambda(1)}\ast.} \\
  K[\Lambda]^{W,\gamma_\xi} \ar[u]^{\subseteq} \ar[r]^{S_1^{-1}} & \Hscr(G,1_{U_0})  }
$$
as well as the inclusion
$$
\Theta(K[\Lambda]^{W,\gamma_\xi}) \subseteq center\; of\;
\Hscr(G,1_{U_1}) \ .
$$
It is known (cf.\ [HKP] Lemma 2.3.1, section 4.6, and Lemma 3.1.1)
that:

-- $\Theta^{norm}(K_q[\Lambda]^W) = center\; of\; \Hscr(G,1_{U_1})
\otimes_K K_q$;

-- $\psi_{\lambda(1)} \ast \Theta^{norm} \circ S^{norm} = id$ on
$\Hscr(G,1_{U_0}) \otimes_K K_q$;

-- $\Theta^{norm} = \iota \circ \Theta^{norm} \circ \iota_\Lambda$
on $K_q[\Lambda]^W$.

The first identity implies the asserted inclusion. We further
deduce that
$$
\eqalign{
 \Theta \circ S_1 & = \iota\circ (\delta^{1/2}\cdot\Theta^{norm})
 \circ\iota_\Lambda\circ (\delta^{1/2}\cdot S^{norm}) \cr
 & = \iota\circ \Theta^{norm}\circ \iota_\Lambda\circ S^{norm} \cr
 & = \Theta^{norm} \circ S^{norm} \qquad\qquad\hbox{on}\
      \Hscr(G,1_{U_0}) \otimes_K K_q }
$$
and hence that
$$
\psi_{\lambda(1)} \ast (\Theta \circ S_1) = id
\qquad\qquad\hbox{on}\ \Hscr(G,1_{U_0}) \ .
$$

\medskip

{\bf 5. Crystalline Galois representations}

\medskip

We go back to the example of the group $G := GL_{d+1}(L)$ which we
have discussed already at the end of section 3. But we now want to
exploit Lemma 2.7. As before we fix a dominant integral weight $\xi
\in X^\ast(T)$ that is given by
$$
\pmatrix{ g_1 &&0 \cr & \ddots \cr 0&& g_{d+1}} \longmapsto
\prod_{i=1}^{d+1} g_i^{a_i}
$$
with $(a_1,\ldots,a_{d+1}) \in \Zdss^{d+1}$. Note that the
dominance means that
$$
a_1 \leq \ldots\leq a_{d+1} \ .
$$
Equally as before we use the coordinates
$$
\matrix{
 T'(K) & \longrightarrow & (K^\times)^{d+1} \hfill\cr
 \hfill \zeta & \longmapsto & (\zeta_1,\ldots,\zeta_{d+1}) \
 \hbox{with}\ \zeta_i := q^{i-1}\pi_L^{a_i}\zeta(\lambda_{\{i\}})
 }
$$
on the dual torus. Some times we view $\zeta$ as the diagonal matrix
in $GL_{d+1}(K)$ with diagonal entries
$(\zeta_1,\ldots,\zeta_{d+1})$. On the other hand, on the root space
we use the coordinates
$$
\matrix{
 V_\Rdss = \Hom(\Lambda,\Rdss) & \longrightarrow & \Rdss^{d+1} \hfill\cr
 \hfill z & \longmapsto & (z_1,\ldots,z_{d+1}) \
 \hbox{with}\ z_i := z(\lambda_{\{i\}}) \ .
 }
$$
In these coordinates we have:

1) The points $\eta_L$ and $\xi_L$ from Example 3 in section 2
correspond to
$$
{[L:\Qdss_p] \over 2}(-d,-(d-2),\ldots,d-2,d)\quad\hbox{and}\quad
(a_1,\ldots,a_{d+1})\ ,
$$
respectively.

2) The map $val : T'(K) \longrightarrow V_\Rdss$ corresponds to
the map
$$
\matrix{
 \hfill (K^\times)^{d+1} & \longrightarrow & \Rdss^{d+1} \hfill\cr
 (\zeta_1,\ldots,\zeta_{d+1}) & \longmapsto &
 (\val_L(\zeta_1),\ldots,\val_L(\zeta_{d+1})) - \xi_L -
 \widetilde{\eta}_L \ . }
$$
where
$$
\widetilde{\eta}_L := [L:\Qdss_p](0,1,\ldots,d) = \eta_L +
{[L:\Qdss_p] \over 2}(d,\ldots,d) \ .
$$

3) On $\Rdss^{d+1}$ the partial order defined by $\Phi^+$ is given
by
$$
(z_1,\ldots,z_{d+1}) \leq (z'_1,\ldots,z'_{d+1})
$$
if and only if
$$
z_{d+1} \leq z'_{d+1}\ ,\ z_d + z_{d+1} \leq z'_d + z'_{d+1}\ ,\
\ldots,\ z_2 +\ldots + z_{d+1} \leq z'_2 + \ldots + z'_{d+1}
$$
and
$$
z_1 + \ldots + z_{d+1} = z'_1 + \ldots + z'_{d+1} \ .
$$

4) The map $z \longmapsto z^{dom}$ corresponds in $\Rdss^{d+1}$ to
the map which rearranges the coordinates in increasing order and
which we will also denote by $(.)^{dom}$.

It now is a straightforward computation to show that Lemma 2.7
amounts to
$$
T'_\xi = \{\zeta \in T' :
(\val_L(\zeta_1),\ldots,\val_L(\zeta_{d+1}))^{dom}
 \leq \xi_L + \widetilde{\eta}_L\} .
$$

For any increasing sequence $\underline{r} = (r_1 \leq\ldots\leq
r_{d+1})$ of real numbers we denote by $\Pscr(\underline{r})$ the
convex polygon in the plane through the points
$$
(0,0),(1,r_1),(2,r_1 + r_2),\ldots,(d+1,r_1+\ldots +r_{d+1}) \ .
$$
We then may reformulate the above description of $T'_\xi$ as
follows.

\medskip

{\bf Lemma 5.1:} {\it $T'_\xi$ is the subdomain of all $\zeta \in
T'$ such that $\Pscr(val(\zeta)^{dom})$ lies above $\Pscr(\xi_L +
\widetilde{\eta}_L)$ and both polygons have the same endpoint.}

\medskip

We recall that a filtered $K$-isocrystal is a triple $\underline{D}
= (D,\varphi, Fil^\cdot D)$ consisting of a finite dimensional
$K$-vector space $D$, a $K$-linear automorphism $\varphi$ of $D$ --
the ``Frobenius'' -- , and an exhaustive and separated decreasing
filtration $Fil^\cdot D$ on $D$ by $K$-subspaces. In the following
we fix the dimension of $D$ to be equal to $d+1$ and, in fact, the
vector space $D$ to be the $d+1$-dimensional standard vector space
$D = K^{d+1}$. We then may think of $\varphi$ as being an element in
the group $G'(K) := GL_{d+1}(K)$. The (filtration) type
$type(\underline{D}) \in \Zdss^{d+1}$ is the sequence
$(b_1,\ldots,b_{d+1})$, written in increasing order, of the break
points $b$ of the filtration $Fil^\cdot D$ each repeated $dim_K\,
gr^b D$ many times. We put
$$
t_H(\underline{D}) := \sum_{b \in \Zdss} b\cdot dim_K\, gr^b\, D \ .
$$
Then $(d+1,t_H(\underline{D}))$ is the endpoint of the polygon
$\Pscr(type(\underline{D}))$. On the other hand we define the
Frobenius type $s(\underline{D})$ of $\underline{D}$ to be the
conjugacy class of the semisimple part of $\varphi$ in $G'(K)$. We
put
$$
t^L_N(\underline{D}) := \val_L ({\rm det}_K(\varphi)) \ .
$$
The filtered $K$-isocrystal $\underline{D}$ is called weakly
$L$-admissible if $t_H(\underline{D}) = t^L_N(\underline{D})$ and
$t_H(\underline{D}') \leq t^L_N(\underline{D}')$ for any filtered
$K$-isocrystal $\underline{D}'$ corresponding to a
$\varphi$-invariant $K$-subspace $D' \subseteq D$ with the induced
filtration.

\medskip

{\bf Proposition 5.2:} {\it Let $\zeta \in T'(K)$ and let $\xi$ be a
dominant integral weight of $G$; then $\zeta \in T'_\xi(K)$ if and
only if there is a weakly $L$-admissible filtered $K$-isocrystal
$\underline{D}$ such that $type(\underline{D}) = \xi_L +
\widetilde{\eta}_L$ and $\zeta \in s(\underline{D})$.}

Proof: Let us first suppose that there exists a filtered
$K$-isocrystal $\underline{D}$ with the asserted properties. Then
$\Pscr(type(\underline{D})) = \Pscr(\xi_L + \widetilde{\eta}_L)$ is
the Hodge polygon of $\underline{D}$ and $\Pscr(val(\zeta)^{dom})$
is its Newton polygon (relative to $\val_L$). By [Fon] Prop.\ 4.3.3
(the additional assumptions imposed there on the field $K$ are
irrelevant at this point) the weak admissibility of $\underline{D}$
implies that its Newton polygon lies above its Hodge polygon with
both having the same endpoint. Lemma 5.1 therefore implies that
$\zeta \in T'_\xi(K)$.

We now assume vice versa that $\zeta \in T'_\xi(K)$. We let
$\varphi_{ss}$ be the semisimple automorphism of the standard vector
space $D$ given by the diagonal matrix with diagonal entries
$(\zeta_1,\ldots,\zeta_{d+1})$. Let
$$
D = D_1 + \ldots + D_m
$$ be the decomposition of $D$ into the eigenspaces of
$\varphi_{ss}$. We now choose the Frobenius $\varphi$ on $D$ in
such a way that $\varphi_{ss}$ is the semisimple part of $\varphi$
and that any $D_j$ is $\varphi$-indecomposable. In this situation
$D$ has only finitely many $\varphi$-invariant subspaces $D'$ and
each of them is of the form
$$
D' = D'_1 + \ldots + D'_m
$$
with $D'_j$ one of the finitely many $\varphi$-invariant subspaces
of $D_j$. By construction the Newton polygon of $(D,\varphi)$ is
equal to $\Pscr(val(\zeta)^{dom})$. To begin with consider any
filtration $Fil^\cdot D$ of type $\xi_L + \widetilde{\eta}_L$ on $D$
and put $\underline{D} := (D,\varphi, Fil^\cdot D)$. The
corresponding Hodge polygon then is $\Pscr(\xi_L +
\widetilde{\eta}_L)$. By Lemma 5.1 the first polygon lies above the
second and both have the same endpoint. The latter already says that
$$
t_H(\underline{D}) = t^L_N(\underline{D}) \ .
$$
It remains to be seen that we can choose the filtration $Fil^\cdot
D$ in such a way that $t_H(\underline{D}') \leq
t^L_N(\underline{D}')$ holds true for any of the above finitely many
$\varphi$-invariant subspaces $D' \subseteq D$. The inequality
between the two polygons which we have does imply that
$$
a_1 + (a_2 +[L:\Qdss_p]) + \ldots + (a_{dim D'} + (dim\, D' -
1)[L:\Qdss_p]) \leq t^L_N(\underline{D}') \ .
$$
Hence it suffices to find the filtration in such a way that we
have
$$
t_H(\underline{D}') \leq a_1 + (a_2 +[L:\Qdss_p]) + \ldots + (a_{dim
D'} + (dim\, D' - 1)[L:\Qdss_p])
$$
for any $D'$. But it is clear that for any filtration (of type
$\xi_L + \widetilde{\eta}_L$) in general position we actually have
$$
t_H(\underline{D}') = a_1 + (a_2 +[L:\Qdss_p]) + \ldots + (a_{dim
D'} + (dim\, D' - 1)[L:\Qdss_p])
$$
for the finitely many $D'$.

\medskip

In order to connect this to Galois representations we have to begin
with a different kind of filtered isocrystal (cf.\ [BM] \S3.1).
First of all we now suppose that $K$ is a finite extension of
$\Qdss_p$ (as always containing $L$). Then a filtered isocrystal
over $L$ with coefficients in $K$ is a triple $\underline{M} =
(M,\phi, Fil^\cdot M_L)$ consisting of a free $L_0 \otimes_{\Qdss_p}
K$-module $M$ of finite rank, a $\sigma$-linear automorphism $\phi$
of $M$ -- the ``Frobenius'' -- , and an exhaustive and separated
decreasing filtration $Fil^\cdot M_L$ on $M_L := L \otimes_{L_0} M$
by $L \otimes_{\Qdss_p} K$-submodules. Here $L_0$ denotes the
maximal unramified subextension of $L$ and $\sigma$ its Frobenius
automorphism. By abuse of notation we also write $\sigma$ for the
automorphism $\sigma \otimes id$ of $L_0 \otimes_{\Qdss_p} K$. We
put
$$
t_H(\underline{M}) := \sum_{b \in \Zdss} b\cdot dim_L\, gr^b\, M_L =
[K:L] \cdot \sum_{b \in \Zdss} b\cdot dim_K\, gr^b\, M_L \ .
$$
The equality is a consequence of the fact that for any finitely
generated $L \otimes_{\Qdss_p} K$-module $M'$ the identity
$$
dim_L M' = [K:L]\cdot dim_K M'
$$
holds true. By semisimplicity this needs to be verified only for a
simple module which must be isomorphic to a field into which $L$ and
$K$ both can be embedded and in which case this identity is obvious.

The number $t_N(\underline{M})$ is defined as
$\val_{\Qdss_p}(\phi(x)/x)$ where $x$ is an arbitrary nonzero
element in the maximal exterior power of $M$ as an $L_0$-vector
space. But we have
$$
\eqalign{
 t_N(\underline{M}) & = \val_{\Qdss_p}(\phi(x)/x)\cr
  & = {1 \over [L_0 : \Qdss_p]} \cdot
 \val_{\Qdss_p}({\rm det}_{L_0}(\phi^{[L_0 : \Qdss_p]})) \cr
 & = {1 \over [L_0 : \Qdss_p]} \cdot
 \val_{\Qdss_p}({\rm Norm}_{K/L_0}({\rm det}_K(\phi^{[L_0 : \Qdss_p]}))) \cr
 & = \val_{\Qdss_p}({\rm Norm}_{K/L_0}({\rm det}_K(\phi))) \cr
 & = [K:L_0] \cdot \val_{\Qdss_p}({\rm det}_K(\phi)) \cr
 & = [K:L] \cdot \val_L({\rm det}_K(\phi)) \ . }
$$
The filtered isocrystal $\underline{M}$ over $L$ with coefficients
in $K$ is called weakly admissible (cf.\ [BM] Prop.\ 3.1.1.5) if
$t_H(\underline{M}) = t_N(\underline{M})$ and $t_H(\underline{M}')
\leq t_N(\underline{M}')$ for any subobject $\underline{M}'$ of
$\underline{M}$ corresponding to a $\phi$-invariant $L_0
\otimes_{\Qdss_p} K$-submodule $M' \subseteq M$ with the induced
filtration on $L \otimes_{L_0} M'$.

By the main result of [CF] there is a natural equivalence of
categories $V \longmapsto D_{cris}(V)$ between the category of
$K$-linear crystalline representations of the absolute Galois group
${\rm Gal}(\overline{L}/L)$ of the field $L$ and the category of
weakly admissible filtered isocrystals over $L$ with coefficients in
$K$. It has the property that
$$
dim_K\, V = rank_{L_0 \otimes_{\Qdss_p} K}\, D_{cris}(V) \ .
$$
To avoid confusion we recall that a $K$-linear Galois representation
is called crystalline if it is crystalline as a $\Qdss_p$-linear
representation. We also recall that the jump indices of the
filtration on $D_{cris}(V)_L$ are called the Hodge-Tate coweights of
the crystalline Galois representation $V$ (they are the negatives of
the Hodge-Tate weights). Moreover, we will say that $V$ is $K$-split
if all eigenvalues of the Frobenius on $D_{cris}(V)$ are contained
in $K$. This is a small technical condition which always can be
achieved by extending the coefficient field $K$. More important is
the following additional condition. We let $\Cdss_p$ denote the
completion of the algebraic closure $\overline{L}$. We may view $V$
as an $L$-vector space through the inclusion $L \subseteq K$.

\medskip

{\bf Definition:} {\it A $K$-linear crystalline representation $V$
of ${\rm Gal}(\overline{L}/L)$ is called special if the kernel of
the natural map $\Cdss_p \otimes_{\Qdss_p} V \dlongrightarrow
\Cdss_p \otimes_L V$ is generated, as a $\Cdss_p$-vector space, by
its ${\rm Gal}(\overline{L}/L)$-invariants (for the diagonal
action).}

\medskip

On the full subcategory of special crystalline Galois
representations we have a simplified form of the above equivalence
of categories. This is well known (see [FR] Remark 0.3). But since
we have not found any details in the literature we include them here
for the convenience of the reader. We will speak of a $K$-isocrystal
and an isocrystal over $L$ with coefficients in $K$, respectively,
if no filtration is prescribed. Suppose that $(M,\phi)$ is an
isocrystal over $L$ with coefficients in $K$. We then have the
$L_0$-isotypic decomposition
$$
M = \oplus_{\tau \in \Delta}\, M_\tau
$$
where $\Delta := {\rm Gal}(L_0/\Qdss_p)$ and where $M_\tau$ is the
$K$-subspace of $M$ on which $L_0$ acts via the embedding $\tau :
L_0 \hookrightarrow K$. One has
$$
\phi(M_\tau) = M_{\tau\sigma^{-1}}
$$
so that $\phi^f$ with $f := |\Delta|$ is an $L_0 \otimes_{\Qdss_p}
K$-linear automorphism of $M$ which respects the above
decomposition. We see that $(M_1,\phi^f|M_1)$ is a $K$-isocrystal
with $dim_K\, M_1 = rank_{L_0 \otimes_{\Qdss_p}\, K} M$.

\medskip

{\bf Lemma 5.3:} {\it The functor
$$
\matrix{
 \hbox{category of isocrystals over $L$} &
 \mathop{\longrightarrow}\limits^{\sim} & \hbox{category of
 $K$-isocrystals}\cr
 \hbox{with coefficients in $K$} & & \cr
 \hfill (M,\phi) & \longmapsto & (M_1,\phi^f|M_1) \hfill }
$$
is an equivalence of categories.}

Proof: Let $\Iscr$ denote the functor in question. To define a
functor $\Jscr$ in the opposite direction let $(D,\varphi)$ be a
$K$-isocrystal. We put $M := L_0 \otimes_{\Qdss_p} D$ and $\phi :=
(\sigma \otimes 1) \circ \phi'$ with
$$
\phi'|M_\tau := \left\{\matrix{\varphi & \hbox{if}\ \tau =
1,\hfill\cr id \hfill & \hbox{otherwise}.}\right.
$$
Here we have used the $K$-linear composed isomorphism
$$
D \longrightarrow L_0 \otimes_{\Qdss_p} D = M
\mathop{\dlongrightarrow}\limits^{pr} M_1
$$
to transport $\varphi$ from $D$ to $M_1$. At the same time it
provides a natural isomorphism $id \simeq \Iscr\circ\Jscr$. The
opposite natural isomorphism $id \simeq \Jscr\circ\Iscr$ is given by
the composed maps
$$
M_{\sigma^i} \mathop{\longrightarrow}\limits^{\phi^i} M_1
\mathop{\longrightarrow}\limits^{\cong} (L_0 \otimes_{\Qdss_p}
M_1)_1 \mathop{\longrightarrow}\limits^{\sigma^{-i} \otimes
\phi^{-f}} (L_0 \otimes_{\Qdss_p} M_1)_{\sigma^i}
$$
for $0 \leq i \leq f-1$.

\medskip

Suppose now that $M_L$ carries a filtration $Fil^\cdot M_L$ making
$\underline{M} := (M,\phi,Fil^\cdot M_L)$ into a filtered isocrystal
over $L$ with coefficients in $K$. Let
$$
M_L = \oplus_\beta\, M_{L,\beta}
$$
where $\beta$ runs over the ${\rm Gal}(\overline{K}/K)$-orbits in
$\Hom_{\Qdss_p}(L,\overline{K})$ be the $L$-isotypic decomposition
of the $L \otimes_{\Qdss_p} K$-module $M_L$. The filtration on $M_L$
induces a filtration $Fil^\cdot M_{L,\beta}$ on each $M_{L,\beta}$
and by the naturality of the decomposition we have
$$
Fil^\cdot M_L = \oplus_\beta\, Fil^\cdot M_{L,\beta} \ .
$$
Moreover, let $\beta_0$ denote the orbit of the inclusion map $L
\subseteq K$. Then $M_{L,\beta_0}$ is the $K$-subspace of $M_L$ on
which $L$ acts through the inclusion $L \subseteq K$. The composite
map
$$
M_1 \mathop{\longrightarrow}\limits^{\subseteq} M \longrightarrow L
\otimes_{L_0} M = M_L \mathop{\dlongrightarrow}\limits^{pr}
M_{L,\beta_0}
$$
is a $K$-linear isomorphism which we may use to transport the
filtration $Fil^\cdot M_{L,\beta_0}$ to a filtration $Fil^\cdot M_1$
on $M_1$. In this way we obtain the filtered $K$-isocrystal
$\underline{D} := (M_1,\phi^f|M_1,Fil^\cdot M_1)$. Obviously the
full original filtration $Fil^\cdot M_L$ can be recovered from
$Fil^\cdot M_1$ if and only if it satisfies
$$
gr^0 M_{L,\beta} = M_{L,\beta} \qquad\hbox{for any}\ \beta \neq
\beta_0 \ .\leqno{(\ast)}
$$
Let us suppose that the condition $(\ast)$ is satisfied. Since
$gr^0$, by definition, does not contribute to the number $t_H(.)$ we
obviously have
$$
t_H(\underline{M}) = [K:L] \cdot t_H(\underline{D}) \ .
$$
On the other hand, using a normal basis of $L_0$ over $\Qdss_p$ as
well as the inverse functor in the proof of Lemma 5.3, we compute
$$
\eqalign{
 t_N(\underline{M}) & = [K:L] \cdot \val_L ({\rm det}_K(\phi))\cr
 & = [K:L] \cdot \val_L ({\rm det}_K((\sigma \otimes 1)\circ
 (\phi^f|M_1 \oplus id_{M_\sigma} \oplus\ldots\oplus id_{M_{\sigma^{f-1}}})))\cr
 & = [K:L] \cdot \val_L ({\rm det}_K(\phi^f|M_1))\cr
 & = [K:L] \cdot t^L_N(\underline{D}) \ . }
$$
With $\underline{M}$ any of its subobjects also satisfies the
condition $(\ast)$. Moreover, by Lemma 5.3, the subobjects of
$\underline{M}$ are in one to one correspondence with the subobjects
of $\underline{D})$. It follows that $\underline{M}$ is weakly
admissible if and only if $\underline{D}$ is weakly $L$-admissible.
Hence we have the induced equivalence of categories
$$
\matrix{
 \hbox{category of weakly admissible} & & \hbox{category of weakly}\cr
 \hbox{filtered isocrystals over $L$ with} &
 \mathop{\longrightarrow}\limits^{\sim} &  \hbox{$L$-admissible filtered}\cr
 \hbox{coefficients in $K$ satisfying $(\ast)$} & & \hbox{$K$-isocrystals.}\cr
  }
$$
Suppose now that $\underline{M} = D_{cris}(V)$ of some $K$-linear
crystalline representation of ${\rm Gal}(\overline{L}/L)$. By the
general theory of crystalline Galois representations we have the
comparison isomorphism
$$
ker(\Cdss_p \otimes_{\Qdss_p} V \dlongrightarrow \Cdss_p \otimes_L
V) \cong \mathop{\oplus}\limits_{i \in \Zdss} \big(\Cdss_p(-i)
\otimes_L (\mathop{\oplus}\limits_{\beta\neq\beta_0} gr^i
M_{L,\beta})\big) \ .
$$
It is Galois equivariant with ${\rm Gal}(\overline{L}/L)$ acting
diagonally on the left and through the first factors on the right.
For the Galois invariants we therefore obtain the formula
$$
ker(\Cdss_p \otimes_{\Qdss_p} V \dlongrightarrow \Cdss_p \otimes_L
V)^{{\rm Gal}(\overline{L}/L)} \cong
\mathop{\oplus}\limits_{\beta\neq\beta_0} gr^0 M_{L,\beta} \ .
$$
It follows that the isocrystal $D_{cris}(V)$ satisfies the condition
$(\ast)$ if and only if the crystalline Galois representation $V$ is
special. Altogether we obtain that the functor $V \longmapsto
D_{cris}(V)_1$ induces an equivalence of categories
$$
\matrix{
 \hbox{category of $K$-linear special} & & \hbox{category of weakly}\cr
 \hbox{crystalline representations} &
 \mathop{\longrightarrow}\limits^{\sim} &  \hbox{$L$-admissible filtered}\cr
 \hbox{of ${\rm Gal}(\overline{L}/L)$} & & \hbox{$K$-isocrystals.}\cr
  }
$$
It satisfies
$$
dim_K\, V = dim_K\, D_{cris}(V)_1 \ .
$$
Finally suppose that $V$ is a
$$
\matrix{ \hbox{$(d+1)$-dimensional $K$-linear $K$-split special
crystalline}\cr
 \hbox{representation of ${\rm Gal}(\overline{L}/L)$ all of
whose Hodge-Tate coweights have}\cr
 \hbox{multiplicity one and increase at
least by $[L:\Qdss_p]$ in each step.} }\leqno{(+)}
$$
Precisely in this situation there is a dominant integral $\xi =
(a_1,\ldots,a_{d+1})$ such that the Hodge-Tate coweights of $V$ are
$\xi_L + \widetilde{\eta}_L$. By Prop.\ 5.2 we find an up to
permutation unique point $\zeta \in T'_\xi(K)$ such that $\zeta \in
s(D_{cris}(V))$. This means we have constructed a surjection
$$
\hbox{set of isomorphism classes of $V$'s with $(+)$}
\dlongrightarrow \bigcup_{\xi\ dominant}^\cdot W\backslash T'_\xi(K)
\ .
$$
Let us again fix a dominant $\xi = (a_1,\ldots,a_{d+1})$ and let
$\rho_\xi$ denote the irreducible rational representation of $G$ of
highest weight $\xi$. By Prop.\ 2.4 and Cor.\ 3.6 we have an
identification
$$
W\backslash T'_\xi(K) \subseteq (W\backslash T'_\xi)(K) \simeq
Sp(\Bscr(G,\rho_\xi|U_0))(K)
$$
where $Sp(\Bscr(G,\rho|U_0))(K)$ the space of $K$-rational points of
the affinoid variety $\Bscr(G,\rho_\xi|U_0)$, i.e., the space of
$K$-valued characters of the Banach-Hecke algebra
$\Bscr(G,\rho_\xi|U_0)$. Our map therefore becomes a map
$$
\matrix{
 \hbox{set of isomorphism classes of}\cr
 \hbox{$(d+1)$-dimensional $K$-linear $K$-split}\cr
 \hbox{special crystalline representations of}
 & \longrightarrow &
  Sp(\Bscr(G,\rho_\xi|U_0))(K) \cr
 \hbox{${\rm Gal}(\overline{L}/L)$ with Hodge-Tate coweights}\cr
 (a_1, a_2 +[L:\Qdss_p], \ldots, a_{d+1} +
 d[L:\Qdss_p])    }
$$
which we write as $V \longmapsto \zeta(V)$. We point out that in
this form our map is canonical in the sense that it does not depend
on the choice of the prime element $\pi_L$: This choice entered into
our normalization of the Satake map $S_\xi$ and into the coordinates
on $T'$ which we used; it is easy to check that the two cancel each
other out. We also note that in the limit with respect to $K$ this
map is surjective.

We finish this section with a speculation in which way the map which
we have constructed above might be an approximation of a true
$p$-adic local Langlands correspondence. We view a point $\zeta \in
Sp(\Bscr(G,\rho_\xi|U_0))(K)$ as a character $\zeta :
\Bscr(G,\rho_\xi|U_0) \longrightarrow K$. Correspondingly we let
$K_\zeta$ denote the one dimensional $K$-vector space on which
$\Bscr(G,\rho_\xi|U_0)$ acts through the character $\zeta$. We may
``specialize'' the ``universal'' Banach
$\Bscr(G,\rho_\xi|U_0)$-module $B_{U_0}^G(\rho_\xi|U_0)$ from
section 1 to $\zeta$ by forming the completed tensor product
$$
B_{\xi,\zeta} := K_\zeta\,
\widehat{\otimes}_{\Bscr(G,\rho_\xi|U_0)}\, B_{U_0}^G(\rho_\xi|U_0)
\ .
$$
By construction the $K$-Banach space $B_{\xi,\zeta}$ still carries a
continuous and isometric (for the quotient norm) action of $G$. A
future $p$-adic local Langlands correspondence should provide us
with a distinguished correspondence (being essentially bijective)
between the fiber of our map in $\zeta$ (i.e., all $V$ of the kind
under consideration such that $\zeta(V) = \zeta$) and the
isomorphism classes of all topologically irreducible ``quotient''
representations of $B_{\xi,\zeta}$. Unfortunately it is not even
clear that the Banach spaces $B_{\xi,\zeta}$ are nonzero.

In order to describe the existing evidence for this picture we first
have to recall how the characters of the Hecke algebra
$\Hscr(G,1_{U_0})$ can be visualized representation theoretically.
Any element $\zeta \in T'(K)$ can be viewed as a character $\zeta :
T \rightarrow \Lambda \rightarrow K^\times$, and correspondingly we
may form the unramified principal series representation
$$
\matrix{
 \Ind_P^G(\zeta)^\infty := & \hbox{space of all locally constant
functions}\ F : G \longrightarrow K\ \hbox{such that}\hfill\cr &
F(gtn) = \zeta(t)^{-1}F(g)\ \hbox{for any}\ g\in G, t\in T, n\in N }
$$
of $G$. The latter is a smooth $G$-representation of finite length.
By the Iwasawa decomposition $G = U_0P$ the subspace of
$U_0$-invariant elements in $\Ind_P^G(\zeta)^\infty$ is one
dimensional so that the action of $\Hscr(G,1_{U_0})$ on it is given
by a character $\omega_\zeta$. On the other hand $\zeta$ defines in
an obvious way a character of the algebra $K[\Lambda]$ which we also
denote by $\zeta$. Using the Satake isomorphism from section 3 one
then has (cf.\ [Ka1] Lemma 2.4(i))
$$
\omega_\zeta = \zeta \circ S_1 = (\zeta\cdot
\pi_L^{-\val_L(\xi(.))}) \circ S_\xi \ .
$$
By [Ka1] Thm.\ 2.7 the ``specialization'' in $\omega_\zeta$
$$
H_{1,\zeta} := K_{\omega_\zeta}\, \otimes_{\Hscr(G,1_{U_0})}\,
\ind_{U_0}^G(1_{U_0}) \ .
$$
of the ``universal'' $\Hscr(G,1_{U_0})$-module
$\ind_{U_0}^G(1_{U_0})$ from section 1 is an admissible smooth
$G$-representation. Since it also is visibly finitely generated it
is, in fact, of finite length. Since $\ind_{U_0}^G(1_{U_0})$ as a
$G$-representation is generated by its $U_0$-fixed vectors the same
must hold true for any of its quotient representations, in
particular for any quotient of $H_{1,\zeta}$. But the subspace of
$U_0$-invariant vectors in $H_{1,\zeta}$ is one dimensional. It
follows that $H_{1,\zeta}$ possesses  a single irreducible quotient
representation $V_{1,\zeta}$ -- the so called spherical
representation for $\zeta$. One has the $G$-equivariant map
$$
\matrix{
 \hfill H_{1,\zeta} & \longrightarrow & \Ind_P^G(\zeta)^\infty \hfill\cr
 1 \otimes f & \longmapsto & f \ast {\bf 1}_\zeta := \sum_{g \in G/U_0} f(g)g({\bf 1}_\zeta) }
$$
where ${\bf 1}_\zeta \in \Ind_P^G(\zeta)^\infty$ denotes the unique
$U_0$-invariant function with value one in $1 \in G$. Hence
$V_\zeta$ can also be viewed as the, up to isomorphism, unique
irreducible constituent of $\Ind_P^G(\zeta)^\infty$ with a nonzero
$U_0$-fixed vector.

Bringing in again the dominant integral weight $\xi$ we have the
$K$-linear isomorphism
$$
\matrix{
 \ind_{U_0}^G(1_{U_0}) \otimes_K \rho_\xi &
 \mathop{\longrightarrow}\limits^{\cong} &
 \ind_{U_0}^G(\rho_\xi|U_0) \hfill\cr
 \hfill f \otimes x & \longmapsto & f_x(g) := f(g)g^{-1}x \ . }
$$
It is $G$-equivariant if, on the left hand side, we let $G$ act
diagonally. On the left, resp.\ right, hand side we also have the
action of the Hecke algebra $\Hscr(G,1_{U_0})$ through the first
factor, resp.\ the action of the Hecke algebra
$\Hscr(G,\rho_\xi|U_0)$. Relative to the isomorphism
$\iota_{\rho_\xi}$ between these two algebras discussed in section 1
the above map is equivariant for these Hecke algebra actions as
well. (Warning: But this map does not respect our norms on both
sides.) By abuse of notation we will use the same symbol to denote
characters of these two Hecke algebras which correspond to each
other under the isomorphism $\iota_{\rho_\xi}$. We obtain an induced
$G$-equivariant isomorphism
$$
H_{1,\zeta} \otimes_K \rho_\xi
\mathop{\longrightarrow}\limits^{\cong} H_{\xi,\zeta} :=
K_{\omega_\zeta}\, \otimes_{\Hscr(G,\rho_\xi|U_0)}\,
\ind_{U_0}^G(\rho_\xi|U_0)
$$
between ``specializations''. Since with $V_{1,\zeta}$ also
$$
V_{\xi,\zeta} := V_{1,\zeta} \otimes_K \rho_\xi
$$
is irreducible as a $G$-representation ([ST1] Prop.\ 3.4) we see
that $V_{\xi,\zeta}$ is the unique irreducible quotient of
$H_{\xi,\zeta}$ and is also the, up to isomorphism, unique
irreducible constituent of $\Ind_P^G(\zeta)^\infty \otimes_K
\rho_\xi$ which as a $U_0$-representation contains $\rho_\xi|U_0$.

Assuming once more that $\zeta \in T'_\xi(K)$ we, of course, have
that
$$
B_{\xi,\zeta} = \hbox{Hausdorff completion of}\ H_{\xi,\zeta}
$$
with respect to the quotient seminorm from
$\ind_{U_0}^G(\rho_\xi|U_0)$. We remark that the unit ball in
$\ind_{U_0}^G(\rho_\xi|U_0)$ and a fortiori its image in
$H_{\xi,\zeta}$ are finitely generated over the group ring $o_K[G]$.
Hence in order to prove that the quotient topology on
$H_{\xi,\zeta}$ is Hausdorff, i.e., that the canonical map
$H_{\xi,\zeta} \longrightarrow B_{\xi,\zeta}$ is injective it
suffices to exhibit some bounded open $G$-invariant $o_K$-submodule
in $H_{\xi,\zeta}$.

\medskip

{\bf Example 1:} Let $G = GL_2(\Qdss_p)$, $\xi = (a_1,a_2)$ a
dominant weight, and $\zeta = (\zeta_1,\zeta_2) \in (K^\times)^2$.
By the discussion at the end of section 2 the defining conditions
for the affinoid domain $T'_\xi$ are
$$
|\zeta_i| \leq |p|^{a_1}\quad\hbox{for}\ i = 1,2\quad\hbox{and}\quad
|\zeta_1\zeta_2| = |p|^{a_1 + a_2 +1} \ .
$$
The complete list of the weakly $\Qdss_p$-admissible filtered
$K$-isocrystals with a Frobenius $\varphi$ whose semisimple part is
given by $\zeta$ is well known (cf.\ [BB] end of section 3.1): Up to
conjugation we may assume that $|\zeta_1| \geq
|\zeta_2|$.\hfill\break
 {\it Case 1:} $|\zeta_1| = |p|^{a_1}$ and $|\zeta_2| = |p|^{a_2 + 1}$;
 then $\varphi$ is semisimple, and there are (up to isomorphism) exactly two weakly
 $\Qdss_p$-admissible filtrations; one corresponds to a
 decomposable and the other to a reducible but indecomposable Galois
 representation.\hfill\break
 {\it Case 2:} $\zeta_1 \neq \zeta_2$ with $|\zeta_i| < |p|^{a_1}$ for $i =
 1,2$; then $\varphi$ is semisimple, and there is (up to isomorphism) exactly one weakly
 $\Qdss_p$-admissible filtration; it corresponds to an irreducible
 Galois representation.\hfill\break
 {\it Case 3:} $\zeta_1 = \zeta_2$ with $|\zeta_i| < |p|^{a_1}$; then $\varphi$
 is not semisimple, and there is (up to isomorphism) exactly one weakly
 $\Qdss_p$-admissible filtration; it corresponds to an irreducible
 Galois representation.\hfill\break
 In particular, the fiber of our above surjection consists of two elements
 in case 1 and of one element in cases 2 and 3.

On the other hand for $|\zeta_1| \geq |\zeta_2|$ the map
$H_{\xi,\zeta} \mathop{\longrightarrow}\limits^{\cong}
\Ind_P^G(\zeta)^\infty \otimes_K \rho_\xi$ always is an isomorphism.
It therefore follows from [BB] Thm.\ 4.3.1 that our $B_{\xi,\zeta}$
coincides in Case 2 with the representation denoted by $\Pi(V)$ in
loc.\ cit. Moreover, still in Case 2, by [BB] Cor.s 5.4.1/2/3 the
representation of $G$ in the Banach space $B_{\xi,\zeta}$ is
topologically irreducible (in particular nonzero) and admissible in
the sense of [ST2] \S3. In Case 3 the same assertions are shown in
[Bre] Thm.\ 1.3.3 under the restriction that $a_2 - a_1 < 2p - 1$
and $a_1 + a_2 \neq -3$ if $p \neq 2$, resp.\ $a_2 - a_1 < 2$ and
$a_1 + a_2 \neq -1$ if $p = 2$.

We mention that in contrast to $B_{\xi,\zeta}$ the representation
$\Ind_P^G(\zeta)^\infty \otimes_K \rho_\xi$ (or equivalently
$\Ind_P^G(\zeta)^\infty$) is irreducible if and only if $\zeta_2
\neq p \zeta_1$. Hence reducibility can only occur for $a_1 = a_2$
in Case 1 and for $a_1 < a_2$ in Case 2.

It was Breuil's fundamental idea that the two dimensional
crystalline Galois representations of ${\rm
Gal}(\overline{\Qdss}_p/\Qdss_p)$ with distinct Hodge-Tate weights
should correspond to the Banach representations $B_{\xi,\zeta}$ of
$GL_2(\Qdss_p)$. Our general speculation therefore should be seen as
an attempt to extend his picture. But we warn the reader that the
case of $GL_2$ is atypical insofar as in general, given a pair
$(\xi,\zeta)$, there will be infinitely many possibilities for a
weakly admissible filtration.

\medskip

{\bf Example 2.} The unit ball $\ind_{U_0}^G(1_{U_0})^0$ in the
normed space $\ind_{U_0}^G(1_{U_0})$ is a module for the unit ball
$\Hscr(G,1_{U_0})^0$ in the Hecke algebra $\Hscr(G,1_{U_0})$ (for
the sup-norm in both cases). For the two groups $G = GL_2(L)$ and $G
= GL_3(L)$ it is known that $\ind_{U_0}^G(1_{U_0})^0$ is free as an
$\Hscr(G,1_{U_0})^0$-module. For $G = GL_2(L)$ this is a rather
elementary computation on the tree and for $G = GL_3(L)$ it is the
main result in [BO] Thm.\ 3.2.4 (see also the paragraph after Thm.\
1.5; we point out that the arguments in this paper actually prove
freeness and not only flatness). Let $\{b_j\}_{j \in \Ndss}$ be a
basis. Then $\{1 \otimes b_j\}_j$ is a basis of $H_{1,\zeta}$ as a
$K$-vector space, and $\sum_j o_K\cdot(1 \otimes b_j)$ is open in
$H_{1,\zeta}$ for the quotient topology provided $\zeta \in
T'_1(K)$. Hence the quotient topology on $H_{1,\zeta}$ is Hausdorff
which means that the natural map $H_{1,\zeta} \longrightarrow
B_{1,\zeta}$ is injective. In particular, $B_{1,\zeta}$ is nonzero.

\medskip

{\bf Example 3:} Let $G = GL_{d+1}(L)$ be general but assume that
$\zeta \in \Hom(\Lambda,o_K^\times) \subseteq T'(K)$. Then, for any
element $F \in \Ind_P^G(\zeta)^\infty$ the function $|F|(g) :=
|F(g)|$ is right $P$-invariant. Since $G/P$ is compact we therefore
may equip $\Ind_P^G(\zeta)^\infty$ with the $G$-invariant norm
$$
\|F\| := \sup_{g \in G} |F|(g) \ .
$$
Moreover, our above map
$$
\ind_{U_0}^G(1_{U_0}) \longrightarrow H_{1,\zeta} \longrightarrow
\Ind_P^G(\zeta)^\infty
$$
then is continuous. Assuming in addition that $\zeta \in T'_\xi(K)$
we obtain by completion a $G$-equivariant continuous $K$-linear map
$$
B_{1,\zeta} \longrightarrow \Ind_P^G(\zeta)^c \ .
$$
The completion $\Ind_P^G(\zeta)^c$ of $\Ind_P^G(\zeta)^\infty$ is
explicitly given by
$$
\matrix{
 \Ind_P^G(\zeta)^c := & \hbox{space of all continuous
functions}\ F : G \longrightarrow K\ \hbox{such that}\hfill\cr &
F(gtn) = \zeta(t)^{-1}F(g)\ \hbox{for any}\ g\in G, t\in T, n\in N }
$$
It is easy to show that $\Ind_P^G(\zeta)^c$ as a representation of
$G$ in a $K$-Banach space is admissible.

\medskip

{\bf Conjecture:} {\it If $\zeta$ is regular then the representation
of $G$ in the Banach space $\Ind_P^G(\zeta)^c$ is topologically
irreducible.}

\medskip

Suppose therefore that $\zeta$ is regular, i.e., is not fixed by any
$1 \neq w \in W$ for the conjugation action of $W$ on $T'$). It is
then well known that:\hfill\break
 -- The smooth $G$-representation $\Ind_P^G(\zeta)^\infty$ is
 irreducible (for example by the Bernstein-Zelevinsky classification).\hfill\break
 -- The above map $H_{1,\zeta}
 \mathop{\longrightarrow}\limits^{\cong} \Ind_P^G(\zeta)^\infty$ is
 an isomorphism ([Ka1] Thm.\ 3.2 and Remark 3.3 or [Dat] Lemma 3.1).
 \hfill\break
 The latter in particular implies that the quotient topology on
$H_{1,\zeta}$ is Hausdorff and that the map $B_{1,\zeta}
\longrightarrow \Ind_P^G(\zeta)^c$ has dense image. In this context
we also mention, without proof, the following result.

\medskip

{\bf Proposition 5.4:} {\it For any two $\zeta_1, \zeta_2 \in
\Hom(\Lambda,o_K^\times)$ the vector space of all $G$-equivariant
continuous linear maps from $\Ind_P^G(\zeta_1)^c$ to
$\Ind_P^G(\zeta_2)^c$ is zero if $\zeta_1 \neq \zeta_2$ and is
$K\cdot id$ if $\zeta_1 = \zeta_2$.}

\medskip

For $G = GL_2(\Qdss_p)$ the above conjecture follows from a
combination of [ST1] \S4 and [ST3] Thm.\ 7.1.

\medskip

{\bf 6. Weakly admissible pairs and functoriality}

\medskip

In the traditional Langlands program the irreducible smooth
representations of a general group $G$ over $L$ are put into
correspondence with continuous homomorphisms from the Galois group
${\rm Gal}(\overline{L}/L)$ (or rather the Weil-Deligne group of
$L$) into the Langlands dual group $G'$ of $G$. In order to do
something in this spirit in our setting it is useful to slightly
change our point of view which we motivate by looking once again at
the $GL_{d+1}$-case. We started from a dominant weight $\xi \in
X^\ast(T)$ and an element $\zeta \in T'(K)$ in the dual torus.
Viewing $\zeta$, by our particular choice of coordinates, as a
diagonal matrix $\zeta_c$ in $G'(K) = GL_{d+1}(K)$ we considered the
$K$-isocrystals $(K^{d+1},\varphi)$ such that $\zeta_c$ lies in the
conjugacy class of the semisimple part of $\varphi$. The weight
$\xi$ was used to prescribe the type of the filtration which would
make these isocrystals into filtered isocrystals. Our basic result
then was that among all these filtered $K$-isocrystals there is at
least one weakly $L$-admissible one if and only if $\zeta \in
T'_\xi(K)$. Now we observe that $\xi$ actually can be used to define
a model filtration on $K^{d+1}$. Quite generally, for any
$K$-rational cocharacter $\nu : \Gdss_m \longrightarrow G'$ we
decompose $K^{d+1}$ into weight spaces
$$
K^{d+1} = \oplus_{i \in \Zdss} (K^{d+1})_i
$$
with respect to $\nu$ and put
$$
Fil^\cdot_\nu\, K^{d+1} := \oplus_{j \geq i} (K^{d+1})_j \ .
$$
Because of $X^\ast(T) = X_\ast(T') \subseteq X_\ast(G')$ this in
particular applies to $\xi\widetilde{\eta}_L$. Of course, the
filtration $Fil^\cdot_{\xi\tilde{\eta}_L}\, K^{d+1}$ has no reason
to be weakly $L$-admissible. But any other filtration of the same
type as $Fil^\cdot_{\xi\tilde{\eta}_L}\, K^{d+1}$ is of the form
$gFil^\cdot_{\xi\tilde{\eta}_L}\, K^{d+1}\break =
Fil^\cdot_{g(\xi\tilde{\eta}_L)}\, K^{d+1}$ for some $g \in G'(K)$.
Hence we may express our basic result also by saying that, given the
pair $(\xi,\zeta)$, there is a pair $(\nu,\varphi) \in X_\ast(G')(K)
\times G'(K)$ such that\hfill\break
 -- $\nu$ lies in the $G'(K)$-orbit of $\xi\widetilde{\eta}_L$,\hfill\break
 -- the semisimple part of $\varphi$ is conjugate to $\zeta_c$ in $G'(K)$,
 and\hfill\break
 -- the filtered $K$-isocrystal $(K^{d+1},\varphi,Fil^\cdot_\nu\,
 K^{d+1})$ is weakly $L$-admissible\hfill\break
 if and only if $\zeta \in T'_\xi(K)$.

Let now $G$ be again a general $L$-split reductive group. We denote
by $G'$ its Langlands dual group which we consider to be defined
over $L$ as well (cf.\ [Bor]). In particular, $T'$ is a maximal
$L$-split torus in $G'$. We view our dominant $\xi \in X^\ast(T) =
X_\ast(T') \subseteq X_\ast(G')$, as above, as a $K$-rational
cocharacter $\xi : \Gdss_m \longrightarrow G'$ and $\zeta \in T'(K)
\subseteq G'(K)$. For a general pair $(\nu,b) \in X_\ast(G')(K)
\times G'(K)$ we introduce some constructions and terminology which
is borrowed from [RZ] Chap.\ 1. Let $REP_K(G')$ denote the category
of $K$-rational representations of $G'$ and let $FIC_K$ denote the
category of filtered $K$-isocrystals. Both are additive tensor
categories. The pair $(\nu,b)$ gives rise to the tensor functor
$$
\matrix{
 I_{(\nu,b)} : REP_K(G') & \longrightarrow & FIC_K \hfill\cr
 \hfill (\rho,E) & \longmapsto &
 (E,\rho(b),Fil^\cdot_{\rho\circ\nu}\, E) \ . }
$$

\medskip

{\bf Definition:} {\it The pair $(\nu,b)$ is called weakly
$L$-admissible if the filtered $K$-isocrystal $I_{(\nu,b)}(\rho,E)$,
for any $(\rho,E)$ in $REP_K(G')$, is weakly $L$-admissible.}

\medskip

Suppose that $(\nu,b)$ is weakly $L$-admissible. Then $I_{(\nu,b)}$
can be viewed as a functor
$$
I_{(\nu,b)} : REP_K(G') \longrightarrow FIC_K^{L-adm}
$$
into the full subcategory $FIC_K^{L-adm}$ of weakly $L$-admissible
filtered $K$-isocrys- tals which, in fact, is a Tannakian category
(the shortest argument for this probably is to observe that for a
Galois representation the property of being special crystalline is
preserved by tensor products and to use the Colmez-Fontaine
equivalence of categories). Moreover, letting $Rep_K^{con}({\rm
Gal}(\overline{L}/L))$ denote the category of finite dimensional
$K$-linear continuous representations of ${\rm Gal}(\overline{L}/L)$
we know from the last section that the inverse of the functor
$D_{cris}(.)_1$ induces a tensor functor between neutral Tannakian
categories
$$
FIC_K^{L-adm} \longrightarrow Rep_K^{con}({\rm Gal}(\overline{L}/L))
\ .
$$
By composing these two functors we therefore obtain a faithful
tensor functor
$$
\Gamma_{(\nu,b)} : REP_K(G') \longrightarrow Rep_K^{con}({\rm
Gal}(\overline{L}/L))
$$
which possibly is no longer compatible with the obvious fiber
functors. This is measured by a $G'$-torsor over $K$ ([DM] Thm.\
3.2). By Steinberg's theorem ([Ste] Thm.\ 1.9) that $H^1(K^{nr},G')
= 0$ over the maximal unramified extension $K^{nr}$ of $K$ this
torsor is trivial over $K^{nr}$. It follows then from the general
formalism of neutral Tannakian categories ([DM] Cor.\ 2.9, Prop.\
1.13) that the functor $\Gamma_{(\nu,b)}$ gives rise to a
$K^{nr}$-homomorphism in the opposite direction between the affine
group schemes of the two categories which is unique up to
conjugation in the target group. For $REP_K(G')$ this affine group
scheme of course is $G'$ ([DM] Prop.\ 2.8). For $Rep_K^{con}({\rm
Gal}(\overline{L}/L))$ we at least have that the $K$-rational points
of this affine group scheme naturally contain the Galois group ${\rm
Gal}(\overline{L}/L)$. Hence by restriction we obtain a continuous
homomorphism of groups
$$
\gamma_{\nu,b} : {\rm Gal}(\overline{L}/L) \longrightarrow
G'(K^{nr})
$$
which is determined by the functor $\Gamma_{(\nu,b)}$ up to
conjugation in $G'(K^{nr})$. So we see that any weakly
$L$-admissible pair $(\nu,b)$ determines an isomorphism class of
``Galois parameters'' $\gamma_{\nu,b}$. We remark that if the
derived group of $G'$ is simply connected Kneser ([Kne]) showed that
$H^1(K,G') = 0$ so that in this case the Galois parameter
$\gamma_{\nu,b}$ already has values in $G'(K)$. Following [RZ] p.\
14 and [Win] one probably can establish an explicit formula for the
cohomology class in $H^1(K,G')$ of the torsor in question.

We indicated already earlier that Langlands functoriality (for
smooth representations) requires to work with the normalized Satake
isomorphism $S^{norm}$. This forces us to assume in this section
that our coefficient field $K$ contains a square root of $q$ and to
pick one once and for all. As a consequence we also have a preferred
square root $\delta^{1/2} \in T'(K)$ of $\delta \in T'(K)$. Being
able to work with the normalized Satake map we do not have to
consider the twisted $W$-action on $K[\Lambda]$. But, of course, we
still have a norm in the picture which depends on $\xi$ and which is
the following. We consider the automorphism of $K$-algebras
$$
\matrix{
 a_\xi : \qquad K[\Lambda] & \longrightarrow & K[\Lambda] \hfill\cr
 \hfill \lambda = \lambda(t) & \longmapsto & \delta^{1/2}(\lambda)
 \pi_L^{\val_L(\xi(t))} \lambda }
$$
which intertwines the conjugation action by $W$ on the source with
the twisted action on the target. Pulling back along $a_\xi$ the
norm $\|\ \|_{\gamma_\xi}$ gives the norm
$$
\|\sum_{\lambda \in \Lambda} c_\lambda \lambda\|_\xi^{norm} :=
\mathop{\sup}\limits_{\lambda = \lambda(t)}
|\delta^{1/2}({^w}\lambda) \pi_L^{\val_L(\xi({^w}t))} c_\lambda|
$$
on $K[\Lambda]$ with $w \in W$ for each $\lambda$ being chosen in
such a way that ${^w}\lambda \in \Lambda^{--}$. Let $K\langle
\Lambda;\xi \rangle$ denote the corresponding Banach algebra
completion of $K[\Lambda]$. It follows from Prop.\ 2.4 that
$K\langle \Lambda;\xi \rangle$ is the affinoid algebra of the
affinoid subdomain $T'_{\xi,norm}$ obtained by pulling back $T'_\xi$
along $a_\xi$. Since $a_\xi$ induces on $T'$ the map $\zeta
\longmapsto \delta^{1/2} \pi_L^{\val_L \circ \xi} \zeta$ we deduce
from Lemma 2.7 that
$$
T'_{\xi,norm}(K) = val^{-1}(V_\Rdss^{\xi,norm})
$$
with
$$
V_\Rdss^{\xi,norm} := \{z \in V_\Rdss : z^{dom} \leq \eta_L + \xi_L
\} \ .
$$
We have the commutative diagram
$$
\xymatrix{
    &  \Bscr(G,\rho_\xi|U_0) \ar[d]   \\
                & \|\ \|_\xi\hbox{-completion of}\ \Hscr(G,1_{U_0})
                 \ar[dl]^{S^{norm}} \ar[dr]_{S_\xi}             \\
 K\langle \Lambda;\xi \rangle^W \ar[rr]^{a_\xi} & &
 K\langle \Lambda;\gamma_\xi \rangle^{W,\gamma_\xi}        }
$$
in which, as a consequence of Lemma 3.3 and Prop.\ 3.5, all maps are
isomorphisms of Banach algebras. In this section we use the left
hand sequence of arrows to identify $\Bscr(G,\rho_\xi|U_0)$ with the
algebra of analytic functions on the affinoid space $W\backslash
T'_{\xi,norm}$. In particular, this identifies $(W\backslash
T'_{\xi,norm})(K)$ with the set of $K$-valued (continuous)
characters of the Banach-Hecke algebra $\Bscr(G,\rho_\xi|U_0)$.

\medskip

{\bf Remark:} {\it Using that $\delta(\lambda_{\{i\}}) =
q^{-d+2(i-1)}$ the statement of Prop.\ 5.2 for the group $G =
GL_{d+1}(L)$ becomes: $\zeta \in T'_{\xi,norm}(K)$ if and only if
there is a weakly $L$-admissible filtered $K$-isocrystal
$\underline{D}$ such that $type(\underline{D}) = \xi_L +
\widetilde{\eta}_L$ and the semisimple part of its Frobenius is
given by the diagonal matrix with entries
$q^{d/2}\zeta(\lambda_{\{i\}})$.}

\medskip

We note that in the case where $\eta_L$ happens to be integral
(i.e., if $d[L:\Qdss_p]$ is even) we can go one step further, can
remove completely normalizations accidental to the group
$GL_{d+1}(L)$, and can restate the above remark equivalently as
follows. We have $\zeta \in T'_{\xi,norm}(K)$ if and only if there
is a weakly $L$-admissible filtered $K$-isocrystal $\underline{D}$
such that $type(\underline{D}) = \xi_L + \eta_L$ and the semisimple
part of its Frobenius is given by the diagonal matrix with entries
$\zeta(\lambda_{\{i\}})$. Passing now to a general $G$ this
unfortunately forces us at present to work under the technical
hypothesis that $\eta_L \in X^\ast(T) = X_\ast(T')$. This, for
example, is the case if $[L:\Qdss_p]$ is even or if the group $G$ is
semisimple and simply connected. To emphasize that $\eta_L$ then
will be considered primarily as a rational cocharacter of $T'$ we
will use multiplicative notation and write $\xi\eta_L$ for the
product of the rational cocharacters $\xi$ and $\eta_L$. In this
setting and for general $G$ the analog of Prop.\ 5.2 is the
following.

\medskip

{\bf Proposition 6.1:} {\it Suppose that $\eta_L$ is integral, let
$\xi \in X^\ast(T)$ be dominant, and let $\zeta \in T'(K)$; then
there exists a weakly $L$-admissible pair $(\nu,b)$ (and hence a
Galois parameter $\gamma_{\nu,b}$) such that $\nu$ lies in the
$G'(K)$-orbit of $\xi\eta_L$ and $b$ has semisimple part $\zeta$ if
and only if $\zeta \in T'_{\xi,norm}(K)$.}

Proof: First let $(\nu,b)$ be a weakly $L$-admissible pair as in the
assertion. Further let $\rho : G' \longrightarrow GL(E)$ be any
$K$-rational representation. We then have the weakly $L$-admissible
filtered $K$-isocrystal $(E,\rho(b), Fil^\cdot_{\rho\circ\nu}\, E)$.
Furthermore $\rho\circ\nu$ is conjugate to $\rho\circ (\xi\eta_L)$
in $GL(E)(K)$ and $\rho(\zeta)$ is the semisimple part of $\rho(b)$.
We fix a $K$-rational Borel subgroup $P_E \subseteq GL(E)$ and a
maximal $K$-split torus $T_E \subseteq P_E$ such that $\rho(\zeta)
\in T_E(K)$. There is a unique $K$-rational cocharacter
$(\rho\circ\nu)^{dom} : \Gdss_m \longrightarrow T_E$ which is
dominant with respect to $P_E$ and which is conjugate to
$\rho\circ\nu$ in $GL(E)(K)$. Then $(\rho\circ\nu)^{dom} =
(\rho\circ (\xi\eta_L))^{dom}$ corresponds to the type of the
filtration $Fil^\cdot_{\rho\circ\nu}\, E$ in the sense of section 5.
As in the first part of the proof of Prop.\ 5.2 we know from [Fon]
Prop.\ 4.3.3 that the weak $L$-admissibility of our filtered
isocrystal implies that the Newton polygon
$\Pscr((\rho(val(\zeta)))^{dom})$ lies above the Hodge polygon
$\Pscr((\rho\circ (\xi\eta_L))^{dom})$ with both having the same
endpoint. But, as discussed before Lemma 5.1, this means that
$$
(\rho(val(\zeta)))^{dom} \leq (\rho\circ (\xi\eta_L))^{dom} \ .
$$
According to [FR] Lemma 2.1 the latter implies that
$$
val(\zeta)^{dom} \leq (\xi\eta_L)^{dom} = \xi\eta_L \quad ,\;
\hbox{i.e., that}\quad \zeta \in T'_{\xi,norm}(K) \ .
$$
For the reverse implication we first recall that, given any pair
$(\nu,b)$ and any $K$-rational representation $\rho : G'
\longrightarrow GL(E)$, the associated filtered $K$-isocrystal
$(E,\rho(b),Fil^\cdot_{\rho\circ\nu}\, E)$ carries the canonical
HN-filtration by subobjects (cf.\ [RZ] Prop.\ 1.4). The latter is
stabilized by a unique parabolic subgroup $P^\rho_{(\nu,b)}
\subseteq GL(E)$. We obviously have
$$
\rho(b) \in P^\rho_{(\nu,b)}(K) \ .
$$
The HN-filtrations, being functorial, equip our functor
$I_{(\nu,b)}$ in fact with the structure of an exact
$\otimes$-filtration in the sense of [Saa] IV.2.1.1. The exactness
is trivial since the category $REP_K(G')$ is semisimple. The
compatibility with the tensor product is a theorem of Faltings and
Totaro (independently). It then follows from [Saa] Prop.\ IV.2.2.5
and Thm.\ IV.2.4 that
$$
P_{(\nu,b)} := \bigcap_\rho \rho^{-1}(P^\rho_{(\nu,b)})
$$
is a $K$-rational parabolic subgroup of $G'$.  Since [Saa] only
considers filtrations indexed by integers this requires the
following additional observation. The category $REP_K(G')$ has a
generator ([Saa] II.4.3.2) and is semisimple. From this one deduces
that the jump indices in the HN-filtrations on all the values of our
functor can be written with a common denominator. Hence all these
HN-filtrations can be reindexed simultaneously in such a way that
they become integral, and [Saa] applies. We emphasize that, denoting
by $\Ddss$ the protorus with character group $\Qdss$, one actually
has a (not unique) $K$-rational homomorphism $\iota_{(\nu,b)} :
\Ddss \longrightarrow G'$ whose weight spaces define the
HN-filtration on the functor $I_{(\nu,b)}$. Its centralizer in $G'$
is a Levi subgroup of $P_{(\nu,b)}$.

Note that we have
$$
b \in P_{(\nu,b)}(K) \ .
$$
After these preliminaries we make our choice of the element $b$.

\medskip

{\bf Lemma 6.2:} {\it There is a regular element $b \in G'(K)$ with
semisimple part $\zeta$.}

Proof: Let $M' \subseteq G'$ denote the connected component of the
centralizer of $\zeta$ in $G'$. We have:\hfill\break
 -- $M'$ is connected reductive ([Ste] 2.7.a); \hfill\break
 -- $M'$ is $K$-split of the same rank as $G'$ (since $T' \subseteq
 M'$); \hfill\break
 -- $\zeta \in T'(K) \subseteq M'(K)$; in fact, $\zeta$ lies in the
 center of $M'$. \hfill\break
 The regular unipotent conjugacy class in $M'$, by its unicity
([Ste] Thm.\ 3.3), is defined over $K$. Since $M'$ is $K$-split it
therefore contains a point $u \in M'(K)$ ([Kot] Thm.\ 4.2). We put
$b := \zeta u \in G'(K)$. The centralizer of $b$ in $G'$ contains
with finite index the centralizer of $u$ in $M'$. Hence $b$ is
regular in $G'$ with semisimple part $\zeta$.

\medskip

We now fix $b \in G'(K)$ to be regular with semisimple part $\zeta$.

\medskip

{\bf Lemma 6.3:} {\it There are only finitely many $K$-rational
parabolic subgroups $Q \subseteq G'$ such that $b \in Q(K)$.}

Proof: Obviously it suffices to prove the corresponding statement
over the algebraic closure $\overline{K}$ of $K$. By [Ste] Thm.\ 1.1
there are only finitely many Borel subgroups $Q_0 \subseteq G'$ such
that $b \in Q_0(\overline{K})$. Let $Q \subseteq G'$ be any
parabolic subgroup with $b \in Q(\overline{K})$. It suffices to find
a Borel subgroup $Q_0 \subseteq Q$ such that $b \in
Q_0(\overline{K})$. Consider the Levi quotient $\overline{M}$ of $Q$
and the image $\overline{b} \in \overline{M}(\overline{K})$ of $b$.
Then $\overline{b}$ is contained in some Borel subgroup
$\overline{Q}_0 \subseteq \overline{M}$ (cf.\ [Hu1] Thm.\ 22.2) and
we can take for $Q_0$ the preimage of $\overline{Q}_0$ in $Q$.

\medskip

It follows that with $\nu$ varying over the $G'(K)$-orbit $\Xi
\subseteq X_\ast(G')$ of $\xi\eta_L$ the family of parabolic
subgroups $P_{(\nu,b)}$ actually is finite. Let $P_1, \ldots, P_m$
denote these finitely many parabolic subgroups and write
$$
\Xi = \Xi_1 \cup\ldots\cup \Xi_m \qquad\hbox{with}\ \Xi_i := \{\nu
\in \Xi : P_{(\nu,b)} = P_i\} \ .
$$
We want to show that $\nu \in \Xi$ can be chosen in such a way that
$P_{(\nu,b)} = G'$. Because then the homomorphism $\iota_{(\nu,b)} :
\Ddss \longrightarrow G'$ factorizes through the center of $G'$.
Since by Schur's lemma the center of $G'$ acts through scalars on
any irreducible $K$-rational representation $\rho$ of $G'$ it
follows that the HN-filtration on the filtered isocrystal
$(E,\rho(b),Fil^\cdot_{\rho\circ\nu}\, E)$ for irreducible $\rho$
has only one step. On the other hand, our assumption that $\zeta \in
T'_{\xi,norm}(K)$ together with [FR] Lemma 2.1 imply that this
filtered isocrystal, for any $\rho$, has HN-slope zero. Hence it is
weakly $L$-admissible, first for irreducible $\rho$ and then by
passing to direct sums also for arbitrary $\rho$. This proves that
the pair $(\nu,b)$ is weakly $L$-admissible.

We argue by contradiction and assume that all $P_1,\ldots,P_m \neq
G'$ are proper parabolic subgroups. By [FR] Lemma 2.2.i we then
find, for any $1 \leq i \leq m$, an irreducible $K$-rational
representation $\rho_i : G' \longrightarrow GL(E_i)$ and a $K$-line
$\ell_i \subseteq E_i$ such that
$$
P_i = \hbox{stabilizer in}\ G'\ \hbox{of}\ \ell_i
$$
(in particular, $\ell_i \neq E_i$). We claim that
$P^{\rho_i}_{(\nu,b)}$, for each $\nu \in \Xi_i$, stabilizes the
line $\ell_i$. To see this we have to recall the actual construction
of $\rho_i$ in loc.\ cit. Fix a maximal $K$-split torus $T_i$ in a
Levi subgroup $M_i$ of $P_i$ and fix a Borel subgroup $T_i \subseteq
B_i \subseteq P_i$. By conjugation we may assume that all the
homomorphisms $\iota_{(\nu,b)}$, for $\nu \in \Xi_i$, factorize
through the center of $M_i$. Recall that $M_i$ then is equal to the
centralizer of $\iota_{(\nu,b)}$ in $G'$. Hence we may view these
$\iota_{(\nu,b)}$ as elements in $X_\ast(T_i) \otimes \Qdss$ which
lie in the interior of the facet defined by $P_i$ (the latter
follows from [Saa] Prop.\ IV.2.2.5.1)). Pick on the other hand a
$B_i$-dominant character $\lambda_i \in X^\ast(T_i)$ which lies in
the interior of the facet corresponding to $P_i$ and let $\rho_i$ be
the rational representation of highest weight $\lambda_i$. Then,
according to [FR], the highest weight space $\ell_i \subseteq E_i$
has the required property that $P_i$ is its stabilizer in $G'$. Let
$\lambda \in X^\ast(T_i)$ be any weight in $E_i$ different from
$\lambda_i$. Then $\lambda_i - \lambda$ is a nonzero linear
combination with nonnegative integral coefficients of $B_i$-simple
roots.

{\it Claim:} $(\lambda_i - \lambda)(\iota_{(\nu,b)}) > 0$

Proof: Let $\{\alpha_j : j \in \Delta\} \subseteq X^\ast(T_i)$ be
the set of $B_i$-simple roots and let $J \subseteq \Delta$ denote
the subset corresponding to $P_i$. The highest weight $\lambda_i$
then satisfies
$$
\lambda_i(\check{\alpha}_j) \left\{\matrix{= 0 & \hbox{if}\ j \in
J,\hfill\cr > 0 \hfill & \hbox{if}\ j \not\in J}\right.
$$
where the $\check{\alpha}_j \in X_\ast(T_i)$ denote the simple
coroots. On the other hand the connected center of $M_i$ is equal to
$(\bigcap_{j \in J} ker(\alpha_j))^\circ$, and we have
$$
\alpha_j(\iota_{(\nu,b)}) \left\{\matrix{= 0 & \hbox{if}\ j \in
J,\hfill\cr
> 0 \hfill & \hbox{if}\ j \not\in J.}\right.
$$
We may write
$$
\lambda_i - \lambda = \sum_{j \in \Delta} c_j \alpha_j
\quad\hbox{with}\ c_j \in \Zdss_{\geq 0} \ .
$$
Hence
$$
(\lambda_i - \lambda)(\iota_{(\nu,b)}) = \sum_{j \not\in J} c_j
\alpha_j(\iota_{(\nu,b)}) \geq 0 \leqno{(+)}
$$
and we have to show that $c_j$, for at least one $j \not\in J$, is
nonzero. Let $\lambda' \in X^\ast(T_i)$ denote the unique dominant
element in the orbit of $\lambda$ under the Weyl group of $T_i$.
Then $\lambda'$ also is a weight occurring in $E_i$ and we have
$$
\lambda_i - \lambda' = \sum_{j \in \Delta} d_j \alpha_j \quad
\hbox{and}\quad \lambda' - \lambda = \sum_{j \in \Delta} e_j
\alpha_j \quad\hbox{with}\ d_j,e_j \in \Zdss_{\geq 0} \ .
$$
In particular, $d_j + e_j = c_j$. Suppose first that $\lambda_i \neq
\lambda'$. Then it suffices to find a $j \not\in J$ such that $d_j >
0$. By the proof of [Hum] 13.4 Lemma B we obtain $\lambda'$ from
$\lambda_i$ by successively subtracting simple roots while remaining
inside the weights occurring in $E_i$ in each step. But because of
$\lambda_i(\check{\alpha}_j) = 0$ if $j \in J$ we know ([Hum] 21.3)
that $\lambda_i - \alpha_j$ cannot be a weight occurring in $E_i$
for any $j \in J$. This means of course that we have to have $d_j >
0$ for some $j \not\in J$. Now assume that $\lambda_i = \lambda'$ so
that $\lambda = {^\sigma}\lambda_i$ for some $\sigma$ in the Weyl
group of $T_i$. According to the proof of [Hum] 10.3 Lemma B we
obtain $\lambda$ from $\lambda_i$ in the following way: Let
$\sigma_j$ be the reflection in the Weyl group corresponding to the
simple root $\alpha_j$. Write $\sigma =
\sigma_{j_1}\ldots\sigma_{j_t}$ in reduced form. Then
$$
\lambda_i - \lambda = \sum_{1 \leq s \leq t} \sigma_{j_{s+1}}\ldots
\sigma_{j_t}(\lambda_i)(\check{\alpha}_{j_s})\alpha_{j_s}
$$
with all coefficients being nonnegative integers. Since the
$\sigma_j$ for $j \in J$ fix $\lambda_i$ we may assume that $j_t
\not\in J$. Then the last term in the above sum is
$\lambda_i(\check{\alpha}_{j_t})\alpha_{j_t}$ whose coefficient is
positive.

This claim means that $\ell_i$ is a full weight space of
$\rho_i\circ\iota_{(\nu,b)}$. But it follows from $(+)$ also that
the weight of $\Ddss$ on $\ell_i$ is maximal with respect to the
natural order on the character group $\Qdss$ of $\Ddss$ among all
weights of $\Ddss$ occurring in $E_i$. Hence $\ell_i$ must be the
bottom step in the HN-filtration of the filtered $K$-isocrystal
$\underline{E}_{i,\nu} :=
(E_i,\rho_i(b),Fil^\cdot_{\rho_i\circ\nu}\, E_i)$ for each $\nu \in
\Xi_i$. As such it carries the structure of a subobject
$\underline{\ell}_{i,\nu} \subseteq \underline{E}_{i,\nu}$. As noted
already, due to $\zeta \in T'_{\xi,norm}$, the HN-slope of
$\underline{E}_{i,\nu}$ is zero. By the fundamental property of the
HN-filtration (cf.\ [RZ] Prop.\ 1.4) the HN-slope of
$\underline{\ell}_{i,\nu}$ then must be strictly positive which
means that
$$
t_H(\underline{\ell}_{i,\nu}) > t^L_N(\underline{\ell}_{i,\nu}) \ .
$$
Suppose that we find an $1 \leq i \leq m$ and a $\nu \in \Xi_i$ such
that $\ell_i$ is transversal to the filtration
$Fil^\cdot_{\rho_i\circ\nu}\, E_i$. Let $(a_1,\ldots,a_r)$, resp.\
$(z_1,\ldots,z_r)$, denote the filtration type (in the sense of
section 5), resp.\ the slopes written in increasing order, of the
corresponding $\underline{E}_{i,\nu}$. The transversality means that
$t_H(\underline{\ell}_{i,\nu}) = a_1$. On the other hand, since
$\ell_i$ is a line we must have $t^L_N(\underline{\ell}_{i,\nu}) =
z_j \geq z_1$. But because of $\zeta \in T'_{\xi,norm}(K)$, once
more [FR] Lemma 2.1, and Lemma 5.1 we have $z_1 \geq a_1$ which
leads to the contradictory inequality
$$
t_H(\underline{\ell}_{i,\nu}) \leq t^L_N(\underline{\ell}_{i,\nu}) \
.
$$
It finally remains to justify our choice of $\nu$. Since the
filtration $Fil^\cdot_{\rho_i\circ\nu}\, E_i$ is well defined for
any $\nu \in \Xi$ (and not only $\nu \in \Xi_i$) it suffices to
establish the existence of some $\nu \in \Xi$ such that
$$
\ell_i\ \hbox{is transversal to}\ Fil^\cdot_{\rho_i\circ\nu}\, E_i \
\hbox{for any}\ 1 \leq i \leq m \ .
$$
Let $F_i \subset E_i$ denote the top step of the filtration
$Fil^\cdot_{\rho_i\circ\xi\eta_L}\, E_i$. We have to find an element
$g \in G'(K)$ such that
$$
\rho_i(g)(\ell_i) \not\subseteq F_i \qquad\hbox{for any}\ 1 \leq i
\leq m \ .
$$
For each individual $i$ the set $U_i := \{g \in G' :
\rho_i(g)(\ell_i) \not\subseteq F_i\}$ is Zariski open in $G'$.
Since $\rho_i$ is irreducible the set $U_i$ is nonempty. The
intersection $U := U_1 \cap\ldots\cap U_m$ therefore still is a
nonempty Zariski open subset of $G'$. But $G'(K)$ is Zariski dense
in $G'$ (cf.\ [Hu1] \S34.4). Hence $U$ must contain a $K$-rational
point $g \in U(K)$. Then the cocharacter $\nu := g^{-1}(\xi\eta_L)$
has the properties which we needed.

\medskip

We summarize that, under the integrality assumption on $\eta_L$, any
$K$-valued character of one of our Banach-Hecke algebras
$\Bscr(G,\rho_\xi|U_0)$ naturally gives rise to a nonempty set of
Galois parameters ${\rm Gal}(\overline{L}/L) \longrightarrow
G'(\overline{K})$. The need to pass to the algebraic closure
$\overline{K}$ comes from two different sources: First the element
$\zeta \in T'_{\xi,norm}$ giving rise to a $K$-valued character of
$\Bscr(G,\rho_\xi|U_0)$ in general is defined only over a finite
extension of $K$; secondly, to make Steinberg's theorem applicable
we had to pass to the maximal unramified extension. In the spirit of
our general speculation at the end of the last section we view this
as an approximation to a general $p$-adic Langlands functoriality
principle.

Without the integrality assumption on $\eta_L$ one can proceed at
least half way as follows. Let us fix, more generally, any natural
number $r \geq 1$. We introduce the category of $r$-filtered
$K$-isocrystals $FIC_{K,r}$ whose objects are triples $\underline{D}
= (D,\varphi, Fil^\cdot D)$ as before only that the filtration
$Fil^\cdot D$ is allowed to be indexed by $r^{-1}\Zdss$ (in
particular, $FIC_K = FIC_{K,1}$). The invariants
$t_H(\underline{D})$ and $t^L_N(\underline{D})$ as well as the
notion of weak $L$-admissibility are defined literally in the same
way leading to the full subcategory $FIC^{L-adm}_{K,r}$ of
$FIC_{K,r}$.

\medskip

{\bf Proposition 6.4:} {\it $FIC^{L-adm}_{K,r}$ is a $K$-linear
neutral Tannakian category.}

Proof: This follows by standard arguments from [Tot].

\medskip

The tensor functor
$$
\matrix{
 I_{(\nu,b)} : REP_K(G') & \longrightarrow & FIC_{K,r} \hfill\cr
 \hfill (\rho,E) & \longmapsto &
 (E,\rho(b),Fil^\cdot_{\rho\circ\nu}\, E) }
$$
makes sense for any pair $(\nu,b) \in (X_\ast(G') \otimes
r^{-1}\Zdss)(K) \times G'(K)$ as does the notion of weak
$L$-admissibility of such a pair. With these generalizations Prop.\
6.1 continues to hold in complete generality (involving $2$-filtered
$K$-isocrystals) with literally the same proof. What is missing at
present is the connection between the categories $FIC^{L-adm}_{K,r}$
and $Rep_K^{con}({\rm Gal}(\overline{L}/L))$. This might involve a
certain extension of the Galois group ${\rm Gal}(\overline{L}/L)$.
We hope to come back to this problem in the future.

\bigskip

{\bf References}

\parindent=23truept

\ref{[BO]} Bellaiche J., Otwinowska A.: Platitude du module
universel pour $GL_3$ en caract\'eristique non banale. Bull.\ SMF
131, 507-525 (2003)

\ref{[BB]} Berger L., Breuil C.: Repr\'esentations cristallines
irr\'eductibles de\break $GL_2(\Qdss_p)$. Preprint 2005

\ref{[Bor]} Borel A.: Automorphic $L$-functions. In Automorphic
Forms, Representations and L-Functions. Proc.\ Symp.\ Pure Math.\ 33
(2), pp. 27-61. American Math. Soc. 1979

\ref{[BGR]} Bosch S., G\"untzer U., Remmert R.: Non-Archimedean
Analysis. Ber- lin-Heidelberg-New York: Springer 1984

\ref{[B-GAL]} Bourbaki, N.: Groupes et alg\`ebres de Lie, Chap. 4-6.
Paris: Masson 1981

\ref{[Bre]} Breuil C.: Invariant $\Lscr$ et s\'erie sp\'eciale
$p$-adique. To appear in Ann.\ Sci.\ ENS

\ref{[BM]} Breuil C., M\'ezard A.: Multiplicit\'es modulaires et
repr\'esentations de $GL_2(\Zdss_p)$ et de ${\rm
Gal}(\overline{\Qdss}_p/\Qdss_p)$ en $\ell = p$. Duke Math.\ J.\
115, 205-310 (2002)

\ref{[BT]} Bruhat F., Tits J.: Groupes r\'eductifs sur un corps
local, I, II. Publ. Math. IHES 41, 5-252 (1972), and 60, 5-184
(1984)

\ref{[Car]} Cartier P.: Representations of $\pfr$-adic groups: a
survey. In Automorphic Forms, Representations and L-Functions.
Proc.\ Symp.\ Pure Math.\ 33 (1), pp. 111-155. American Math. Soc.
1979

\ref{[CF]} Colmez P., Fontaine J.-M.: Construction des
repr\'esentations semi-stable. Invent.\ math.\ 140, 1-43 (2000)

\ref{[Dat]} Dat J.-F.: Caract\`eres \`a valeurs dans le centre de
Bernstein. J.\ reine angew.\ Math.\ 508, 61-83 (1999)

\ref{[Del]} Deligne P.: Formes modulaires et representations de
$GL(2)$. In Modular Functions of One Variable II (Eds.\ Deligne,
Kuyk). Lecture Notes in Math.\ 349, pp.\ 55-105.
Berlin-Heidelberg-New York: Springer 1973

\ref{[DM]} Deligne P., Milne J.\ S.: Tannakian categories. In Hodge
Cycles, Motives, and Shimura Varieties (Eds.\ Deligne, Milne, Ogus,
Shih). Lecture Notes in Math.\ 900, pp.\ 101-228.
Berlin-Heidelberg-New York: Springer 1982

\ref{[Em1]} Emerton M.: Jacquet modules of locally analytic
representations of $p$-adic reductive groups I: Definitions and
first properties. To appear in Ann.\ Sci.\ ENS

\ref{[Em2]} Emerton M.: $p$-adic $L$-functions and unitary
completions of representations of $p$-adic reductive groups. To
appear in Duke Math.\ J.

\ref{[Fon]} Fontaine J.-M.: Modules galoisiens, modules filtr\'es et
anneaux Bar- sotti-Tate. Ast\'erisque 65, 3-80 (1979)

\ref{[FR]} Fontaine J.-M., Rapoport M.: Existence de filtrations
admissibles sur des isocristaux. Bull.\ SMF 133, 73-86 (2005)

\ref{[FvP]} Fresnel J., van der Put M.: Rigid Analytic Geometry
and Its Applications. Boston-Basel-Berlin: Birkh\"auser 2004

\ref{[Hai]} Haines T.: The combinatorics of Bernstein functions.
Trans.\ AMS 353, 1251-1278 (2001)

\ref{[HKP]} Haines T., Kottwitz R., Prasad A.: Iwahori-Hecke
algebras. Preprint

\ref{[Gro]} Gross B.H.: On the Satake isomorphism. In Galois
Representations in Arithmetic Algebraic Geometry (Eds.\ Scholl,
Taylor), London Math.\ Soc.\ Lect.\ Notes 254, pp.\ 223-237.
Cambridge Univ.\ Press 1998

\ref{[Hum]} Humphreys J.E.: Introduction to Lie Algebras and
Representation Theory. Berlin-Heidelberg-New York: Springer 1972

\ref{[Hu1]} Humphreys J.E.: Linear Algebraic Groups.
Berlin-Heidelberg-New York: Springer 1987

\ref{[Hu2]} Humphreys J.E.: Reflection groups and Coxeter Groups.
Cambridge Univ.\ Press 1990

\ref{[Jan]} Jantzen J.C.: Representations of Algebraic Groups.
Orlando: Academic Press 1987

\ref{[Ka1]} Kato S.: On Eigenspaces of the Hecke Algebra with
Respect to a Good Maximal Compact Subgroup of a $p$-Adic Reductive
Group. Math.\ Ann.\ 257, 1-7 (1981)

\ref{[Ka2]} Kato S.: Spherical Functions and a $q$-Analogue of
Kostant's Weight Multiplicity Formula. Invent.\ math.\ 66, 461-468
(1982)

\ref{[KKMS]} Kempf G., Knudsen F., Mumford D., Saint-Donat B.:
Toroidal Embeddings I. Lecture Notes in Math.\ 339.
Berlin-Heidelberg-New York: Springer 1973

\ref{[Kne]} Kneser M.: Galois-Kohomologie halbeinfacher
algebraischer Gruppen \"uber $\pfr$-adischen K\"orpern. I. Math.\
Z.\ 88, 40-47 (1965)

\ref{[Kot]} Kottwitz R.: Rational conjugacy classes in reductive
groups. Duke Math.\ J.\ 49, 785-806 (1982)

\ref{[Kut]} Kutzko P.: Mackey's Theorem for non-unitary
representations. Proc. AMS 64, 173-175 (1977)

\ref{[Mac]} Macdonald I.G.: Spherical Functions on a Group of
$p$-Adic Type. Ramanujan Institute Publ.\ 2 (1971)

\ref{[RZ]} Rapoport M., Zink T.: Period Spaces for $p$-divisible
Groups. Annals Math.\ Studies 141. Princeton Univ.\ Press 1996

\ref{[Saa]} Saavedra Rivano N.: Cat\'egories Tannakiennes. Lecture
Notes in Math.\ 265. Berlin-Heidelberg-New York: Springer 1972


\ref{[ST1]} Schneider P., Teitelbaum J.: $U(\gfr)$-finite locally
analytic representations. Representation Theory 5, 111-128 (2001)

\ref{[ST2]} Schneider P., Teitelbaum J.: Banach space
representations and Iwasawa theory. Israel J.\ Math.\ 127, 359-380
(2002)

\ref{[ST3]} Schneider P., Teitelbaum J.: Algebras of $p$-adic
distributions and admissible representations. Invent.\ math.\ 153,
145-196 (2003)

\ref{[Ste]} Steinberg R.: Regular elements of semi-simple algebraic
groups. Publ.\ IHES 25, 49-80 (1965)

\ref{[Tot]} Totaro B.: Tensor products in $p$-adic Hodge theory.
Duke Math.\ J.\ 83, 79-104 (1996)

\ref{[Vig]} Vigneras M.-F.: Alg\`ebres de Hecke affines
g\'en\'eriques. Preprint 2004

\ref{[Win]} Wintenberger J.-P.: Propri\'et\'es du groupe Tannakien
des structures de Hodge $p$-adiques et torseur entre cohomologies
cristalline et \'etale. Ann.\ Inst.\ Fourier 47, 1289-1334 (1997)

\bigskip
\vfill\eject

\parindent=0pt

Peter Schneider\hfill\break Mathematisches Institut\hfill\break
Westf\"alische Wilhelms-Universit\"at M\"unster\hfill\break
Einsteinstr. 62\hfill\break D-48149 M\"unster, Germany\hfill\break
pschnei@math.uni-muenster.de\hfill\break
http://www.uni-muenster.de/math/u/schneider\hfill

\noindent Jeremy Teitelbaum\hfill\break Department of Mathematics,
Statistics, and Computer Science (M/C 249)\hfill\break University of
Illinois at Chicago\hfill\break 851 S. Morgan St.\hfill\break
Chicago, IL 60607, USA\hfill\break jeremy@uic.edu\hfill\break
http://www.math.uic.edu/$\sim$jeremy\hfill

\end